\documentclass[final]{siamltex}

\usepackage{amsmath}
\usepackage{amsfonts}
\usepackage{color}
\usepackage{graphicx}
\usepackage{subfigure}
\def\beq{\begin{equation}}
\def\eeq{\end{equation}}
\def\esplit{\end{split}}
\def\beqalign{\begin{array}{rl}}
\def\eeqalign{\end{array}}

\def\Abold{\mathbf{A}}
\def\Bbold{\mathbf{B}}
\def\Cbold{\mathbf{C}}

\def\Fbold{\mathbf{F}}

\def\Kbold{\mathbf{K}}

\def\Mbold{\mathbf{M}}

\def\Pbold{\mathbf{P}}

\def\Rbold{\mathbf{R}}
\def\Sbold{\mathbf{S}}

\def\Vbold{\mathbf{V}}
\def\Wbold{\mathbf{W}}

\def\Zbold{\mathbf{Z}}

\def\bbold{\mathbf{b}}
\def\cbold{\mathbf{c}}
\def\dbold{\mathbf{d}}

\def\fbold{\mathbf{f}}

\def\qbold{\mathbf{q}}

\def\ubold{\mathbf{u}}

\def\xbold{\mathbf{x}}
\def\ybold{\mathbf{y}}

\def\lambdabold{\boldsymbol{\lambda}}

\def\zerobold{{\bf 0}}
\def\0{\mathbf{\0}}

\newcommand{\bmat}[1]{\begin{bmatrix}#1\end{bmatrix}} 
\newcommand{\pmat}[1]{\begin{pmatrix}#1\end{pmatrix}} 

\title{A practical factorization of a Schur complement for PDE-constrained Distributed Optimal Control\thanks{The first and second authors acknowledge partial support by the Army Research Laboratory through the Army High Performance Computing Research Center under Cooperative Agreement W911NF-07-2-0027. The third and fourth authors acknowledge partial support by the ONR grant N000141110067.}}

\author{Youngsoo Choi\thanks{Aeronautics and Astronautics, Stanford University, Stanford, CA, USA}\and 
				Charbel Farhat\thanks{Aeronautics and Astronautics, Mechanical Engineering, Stanford University, Stanford, CA, USA} \and
				Walter Murray\thanks{Management Science and Engineering, Stanford University, Stanford, CA, USA} \and 
        Michael Saunders\thanks{Management Science and Engineering, Stanford University, Stanford, CA, USA}}

\begin{document}

\maketitle

\begin{abstract}
A distributed optimal control problem with the constraint of a linear elliptic partial differential equation is considered. 
A necessary optimality condition for this problem forms a saddle point system, the efficient and accurate solution 
of which is crucial.
A new factorization of the Schur complement for such a system is proposed and its characteristics discussed.
The factorization introduces two complex factors that are complex conjugate to each other. 
The proposed solution methodology involves the application of a parallel 
linear domain decomposition solver---FETI-DPH---for the solution of the subproblems with the complex factors.
Numerical properties of FETI-DPH in this context are demonstrated,
including numerical and parallel scalability and regularization dependence.
The new factorization can be used to solve Schur complement systems arising 
in both range-space and full-space formulations. In both cases, numerical results indicate
that the complex factorization is promising.
\end{abstract}



\pagestyle{myheadings}
\thispagestyle{plain}
\markboth{Y. CHOI, C. FARHAT, W. MURRAY, AND M. SAUNDERS}{A PRACTICAL FACTORIZATION OF THE SCHUR COMPLEMENT}

	\section{Introduction}
    Numerical methods for solving partial differential equations (PDEs) have broad applications
    in the simulation of complicated physical models, the prediction of their response, and design.
    An important and practical subset of these applications---particularly for design---involves 
    the use of a mathematical optimization technique in which the PDE takes the role of a constraint equation.
		Two methodologies for PDE-constrained optimization problems 
    are SAND (simultaneous analysis and design) and NAND (nested analysis and design) \cite{Biegler2003large}. 
    NAND uses PDEs to express decision variables as an implicit function 
    of state variables and does not include state variables as optimization variables. Thus, the size of
		the optimization problem is not typically large. On the other hand, the SAND approach takes both decision variables
    and state variables as optimization variables and considers PDEs to be equality constraints. Consequently,
    the size of the system in the SAND approach is generally much larger.
    The NAND approach has traditionally been the method of choice for physics-based applications not only because 
    SAND requires the solution of a large-scale system of equations
    but also because NAND conveniently permits the direct use of existing solvers for both optimization 
    and PDE simulation as a black box. However, NAND suffers from the fact that
    many PDE simulations are required for function evaluations---typically the most expensive part of this approach.
    Due to the continuous increase in computational power (e.g., speed of processing, capacity of memory, high performance computing) 
    accompanied by the development of robust and versatile numerical algorithms (e.g., parallel algorithms),
    the SAND approach has received increasing attention from researchers in recent years \cite{Biegler2003sqp, Biros2005I, Choi2012, Prudencio2006},
    and the present study continues this line of research.

    The SAND approach to PDE-constrained optimization takes the form
    \begin{equation}
    \label{eq:pdeconstrainedoptimization}
    \begin{aligned}
      & \underset{y,u}{\text{minimize}}
      & & F(y,u) 
      \\ & \text{subject to}
      & & C(y,u) = 0,
      \end{aligned}
    \end{equation}
    where $C(y,u) = 0$ is the time-independent PDE constraint, $y$ is the vector of state variables, 
		and $u$ is the vector of decision variables.
		State variables by definition are the unknown variables in the forward PDE problem. For example, state variables comprise
		temperature in heat conduction problems and displacements in elastostatics.
    For the class of PDE-constrained optimization known as \emph{optimal control} the decision variables $u$ are referred to as control variables, 
    whereas for \emph{optimal design} or \emph{shape optimization} problems the decision variables $u$ are called design variables. 
    The decision variables may also be a set of parameters describing the material properties or the system of dynamics
    for some inverse problems. All three types of PDE-constrained optimization problems share a similar
		structure of their linear or linearized system of equations known as a KKT system, after the Karush-Kuhn-Tucker optimality conditions, or a saddle point system. 
		
		In what follows, a robust and versatile numerical method for solving a distributed optimal control problem is considered.
		In particular, heat conduction and elastostatic PDE-constrained problems are studied, although
		the method developed here can also be applied or extended easily to other types of PDE-constrained optimization problems.
		From various possible objective functions that may be used to formulate a distributed PDE-constrained optimal control problem, 
    the one considered is that in which a target state is assumed to be given. Thus, the aim is to find a state that
		is close to the prescribed target and a control that realizes that particular state. 
		For example, the static thermal conduction optimal control problem with a target temperature distribution $\bar{y}$ is formulated as
		\begin{equation}
		\label{eq:continuousoptimalheatcontrol}
		\begin{aligned}
		\underset{y,u}{\text{minimize}} & & F(y,u) &:= \frac{1}{2} \int_\Omega(y - \bar{y})^2 dx + \frac{\phi}{2} \int_\Omega u^2 dx
		\\ \text{subject to} & & -\nabla^2y &= u \text{ on } \Omega 
		\\ & & y &= y_c   \text{ on } \Gamma_g.
		\end{aligned}
		\end{equation}
    The solution of this problem is considered in Section~\ref{sec:thermalproblem}.
		Variables $y$ and $u$ are the temperature state and heat source control, respectively, while
		$y_c$ is the prescribed boundary conditions and $\phi$ is a regularization parameter.
    The domain of interest is denoted as $\Omega$ and a Dirichlet boundary condition is imposed on the boundary $\Gamma_g$.
    Note that a unit conductivity is assumed.

		Two alternative approaches to PDE-constrained optimization problems such as (\ref{eq:continuousoptimalheatcontrol}) are
		\textit{optimize-and-discretize} and \textit{discretize-and-optimize}.
		The latter approach, in which one first discretizes both the objective function and constraints 
    and then obtains the discretized optimality condition, is typically used for PDE-constrained optimal control 
    \cite{Biros2005I,Dollar2010,Pearson2012,Pearson2012siam,Rees2010,Simoncini2012,Draganescu2012} and is followed here.

		Solving the saddle point system representing the discretized optimality condition efficiently is crucial 
    to the competitiveness of the SAND approach, in comparison with NAND.
		There are two methodologies for solving a saddle point system.

    \begin{itemize}
      \item	In reduced-space methods, one attempts to reduce the size of a saddle point system by eliminating some variables
		and solving a smaller system \cite{Thorne2009,Thorne2009part2,Simoncini2012,Biegler2003sqp}.
		The range-space method and the null-space method are two popular reduced-space methods.
		In the range-space method, one solves for the dual variables first using the corresponding Schur complement and subsequently updates the primal variables, 
		whereas	the null-space method subdivides the variables algebraically into null-space variables and range-space variables using
		null-space and range-space bases of the constraint Jacobian. 
		  \item In full-space methods, one solves for the primal and dual variables simultaneously.
		The resultant system of equations becomes a sparse saddle point system and
    iterative methods are the only practical choice for large problems. 
    Various efficient preconditioners for saddle point systems have been developed 
		\cite{Murphy2000, Benzi2005, Gill1992, Forsgren2007, Rees2010, Pearson2012, Dollar2006, Dollar2010, 
		Keller2000, Biros2005I, Stoll2008}. However, there is still motivation for further research in this area. 
		For example, most preconditioners are not robust when $\phi$ is small.
		Recently, Pearson and Wathen have developed a new approximation of the Schur complement and used it to facilitate 
		a regularization-robust preconditioner for a particular distributed optimal control problem \cite{Pearson2012}. 
		Pearson et al.\ have further extended the usage of the aforementioned approximation to a broader range of optimal control problems \cite{Pearson2012siam}.
    \end{itemize}

		The factorization of the Schur complement presented in this paper has a similar form to Pearson and Wathen's approximation, 
    and furthermore is applied to the same particular distributed optimal control problem \cite{Pearson2012}.
		However the approach taken here differs from that of Pearson and Wathen in two ways. 
    First and foremost, because the factorization is exact by its nature, the range-space method can be adopted and one efficient 
		solve with the Schur complement results in a solution to the distributed optimal control problem. 
    Second, a scalable domain decomposition based parallel linear solver FETI-DPH \cite{Farhat2005} is used 
    to solve each of the subproblems arising in the application of the proposed factorization. 
    In contrast, Pearson and Wathen suggest a multigrid method.

    The outline of the paper is as follows:
		Section~\ref{sec:distributedcontrolproblems} describes discretization of the distributed optimal control problem 
		\eqref{eq:continuousoptimalheatcontrol} and the corresponding optimality condition, which is a saddle point system.
		Section~\ref{sec:rangefullspacemethod} outlines two methods for solving the saddle point system: 
		the range-space method and the full-space method. Additionally, this section reviews existing preconditioners
		for the full-space method related to the Schur complement in the range-space method.
		Section~\ref{sec:exactrepresentation} introduces a practical factorization of the Schur complement and
		explains how it can be used in both the range-space method and the full-space method.
		Section~\ref{sec:FETI} describes the Finite Element Tearing and Interconnecting (FETI) method and 
    its so-called \emph {dual-primal} variants, FETI-DP and FETI-DPH.
		Section~\ref{sec:numericalresults} presents numerical results that illustrate the scalability and efficiency
		of the method applied to a selection of distributed optimal control problems. 

	\section{Distributed control problems}
	\label{sec:distributedcontrolproblems}
    In this section, the finite element discretization of a distributed optimal control problem is introduced and
    the corresponding optimality condition is presented. For brevity, we choose to
    show the discretization of a distributed optimal control problem with a constraint of the Poisson equation
    (i.e., Eq.~\eqref{eq:continuousoptimalheatcontrol}). However, distributed optimal control problems
    with other types of PDE constraints (e.g., linear elasticity) will have 
    the same discretized formulation. The finite element discretization of \eqref{eq:continuousoptimalheatcontrol} gives
		\begin{equation}
		\label{eq:discretizedoptimalheatcontrol}
		\begin{aligned}
			&\underset{\ybold,\ubold}{\text{minimize}} & & F(\ybold,\ubold) := \frac{1}{2} \|\ybold - \bar{\ybold}\|_\Vbold^2 + \frac{\phi}{2} \|\ubold\|_\Vbold^2
			\\ & \text{subject to} & & \Kbold\ybold + \Kbold_{c}\ybold_c = \Vbold\ubold,
		\end{aligned}
		\end{equation}
		where $\Kbold$, $\Vbold\in \mathbb{R}^{n \times n}$, and $\Kbold_c\in \mathbb{R}^{n \times m}$ are the stiffness matrix, 
		volume matrix, and constrained stiffness matrix, respectively.
		Vector valued quantities $\bar{\ybold}$, $\ybold$, $\ubold \in \mathbb{R}^n$, and 
    $\ybold_c\in \mathbb{R}^m$ are the discretized versions of the target state, 
		state, control, and prescribed boundary conditions, respectively.
    All the discrete variables are denoted with bold fonts for the remainder of the paper.
		The dimensions $n$ and $m$ are the number of unconstrained degrees of freedom (i.e.,\ the size of $\ybold$) 
		and the number of constrained degrees of freedom, respectively.
		It is assumed that both $\Kbold$ and $\Vbold$ are symmetric positive definite (SPD) matrices, which is the case for the thermal problem above.
		The energy norm $\|\cdot\|_V$ is defined as $\|\qbold\|_\Vbold = \sqrt{\qbold^T\Vbold\qbold}$.
		The details of the finite element discretization procedure can be found in \cite{Choi2012}.
		Problem \eqref{eq:discretizedoptimalheatcontrol} is a convex quadratic program and
		its solution is a saddle point of the Lagrangian,
		\begin{equation}
		\label{eq:lagrangian}
			L(\ybold,\ubold,\lambdabold) = \frac{1}{2} \|\ybold - \bar{\ybold}\|_\Vbold^2 + 
                                 \frac{\phi}{2} \|\ubold\|_\Vbold^2 + \lambdabold^T(\Kbold\ybold + \Kbold_{c}\ybold_c - \Vbold\ubold).
		\end{equation}
		A saddle point $(\ybold^*,\ubold^*,\lambdabold^*)$ must satisfy the following optimality condition:	
		\begin{equation}
		\label{eq:linearstaticforcecontrolKKT}
			\bmat{\Vbold    & \zerobold   &  \Kbold
			\\    \zerobold & \phi \Vbold & -\Vbold
			\\    \Kbold    & -\Vbold     &  \zerobold   }
			\pmat{\ybold^* \\ \ubold^* \\ \lambdabold^*} =
			\pmat{\Vbold\bar{\ybold} \\ \zerobold \\ -\Kbold_c \ybold_c}.
		\end{equation}
		This equation is of the form $\Abold\xbold = \bbold$, where $\Abold$ is a symmetric indefinite matrix. 
    Section \ref{sec:rangefullspacemethod} explains two ways of
		solving the saddle point system in \eqref{eq:linearstaticforcecontrolKKT}: the range-space method and full-space method.

	\section{Range-space and full-space methods}
	\label{sec:rangefullspacemethod}

		\subsection{Range-space method}
		\label{subsec:rangespacemethod}
      The range-space method  was introduced for large-scale inequality-constrained convex quadratic programming problems 
      with small active sets in order to overcome the disadvantage of the null-space method:  
      the increasingly high dimension of the null space basis
      as a solution is approached \cite{gill1982rangespace, gill1984weighted}.
      For equality-constrained convex quadratic programming, the range-space method is equivalent to
      the Schur complement method \cite{Thorne2009} where the corresponding dual problem is solved first.
			The range-space method is a reduced-space method in a sense that
			the size of the saddle point system in \eqref{eq:linearstaticforcecontrolKKT} is reduced by eliminating
			$\ybold^*$ and $\ubold^*$ and solving for $\lambdabold^*$ first. The dual variable $\lambdabold^*$ is obtained by solving 
			\begin{equation}
			\label{eq:schurcomplementsystem}
				\Sbold \lambdabold^* = \Kbold_{c}\ybold_c + \Kbold\bar{\ybold},
			\end{equation}
			where $\Sbold = \Kbold\Vbold^{-1}\Kbold+\frac{1}{\phi}\Vbold$ is known as the negative Schur complement on the dual variables.
			Then, $\ybold^*$ and $\ubold^*$ are computed from
			\begin{equation}
			\label{eq:updateyandu}
				\ybold^* = \bar{\ybold}-\Vbold^{-1}\Kbold\lambdabold^* \qquad \text{ and } \qquad
				\ubold^* = \frac{1}{\phi}\lambdabold^*.
			\end{equation}
			The most expensive and crucial step in the range-space method is to solve with $\Sbold$ in \eqref{eq:schurcomplementsystem}.	
      The eigenvalues (and consequently, the condition number) of $\Sbold$ depend on 
      those of $\Kbold$ and $\Vbold$ and the value of $\phi$ \cite{Thorne2009, Choi2012}.
      For small values of $\phi$, the condition number of $\Sbold$ is affected most by the eigenvalues of $\Vbold$, 
      whose condition number is bounded above by a constant, $C$,
      for $P_q$ or $Q_q$ (the $q^{th}$ order triangular or quadrilateral finite elements in 2D and 
      tetrahedral or brick elements in 3D) if a set of grids is quasi-uniform (see Eq.~(1.116) and Eq.~(1.117) in \cite{Elman2005}).
      The constant $C$ depends on the order of approximation $q$ but not on the mesh size $h$.
      Thus for small values of $\phi$, even if the grid is refined or the problem size is increased, the condition number of $\Sbold$ is bounded.
      For large values of $\phi$, the condition number of $\Sbold$ (denoted by $\kappa(\Sbold)$) depends on $\kappa(\Kbold\Vbold^{-1}\Kbold)$,
      which is proportional to $\kappa(\Kbold)^2$. For a second-order elliptic operator, 
      one can prove that $\kappa(\Kbold) < ch^{-2}$ if $P_q$ or $Q_q$ is used on a quasi-uniform discretization of a domain 
      (see Eq.~(1.119) and Eq.~(1.121) in \cite{Elman2005}).
      For a higher-order elliptic operator, the dependence on $h$ of $\kappa(\Kbold)$ increases (e.g., the shell or beam elements).
      This is unfortunate because $\kappa(\Kbold)$ is likely to increase as the mesh size $h$ decreases or the problem size increases.
      One remedy is to apply iterative refinement with a high precision data type.
      If the solution is inaccurate even after iterative refinement due to the ill-conditioning of $\Sbold$,
      an alternative way of solving \eqref{eq:linearstaticforcecontrolKKT} 
			that avoids an exact solve with $\Sbold$ is to use the full-space method, which is explained in the next section.

      An additional challenge in applying the range-space method for solving \eqref{eq:linearstaticforcecontrolKKT} resides in
      solving with $\Sbold$ itself. Because $\Sbold$ is a sum of two matrices ($\Kbold\Vbold^{-1}\Kbold$ and $\frac{1}{\phi}\Vbold$), the first of which is a product 
      requiring two matrix-matrix multiplications to evaluate, it is not straightforward to come up with an efficient way of solving 
      with $\Sbold$. In Section~\ref{sec:exactrepresentation}, a practical factorization of $\Sbold$ is introduced and 
      the factorization does not require matrix-matrix products.
      The factorization introduced in Section~\ref{sec:exactrepresentation} introduces complex symmetric matrices (not Hermitian matrices)
      in its factors. A list of efficient solvers for complex symmetric matrices includes CG-type methods \cite{Freund1992}, 
      CS-MINRES-QLP \cite{choi2013minimal}, GMRES \cite{Saad1986, Day2001}, 
      FETI-DPH \cite{Farhat2005}, and multi-grid methods \cite{Lahye2000, reitzinger2003algebraic}. 
      FETI-DPH is explained in Section \ref{sec:FETI} and used in numerical experiments in Section \ref{sec:numericalresults}. 

		\subsection{Full-space Method}
		\label{sebsec:fullspacemethod}	
			The full-space method attempts to solve for all the variables of 
			the saddle point system in \eqref{eq:linearstaticforcecontrolKKT} simultaneously.
			As the discretization is refined, the size of the system becomes large and
			an iterative method is often the only available method. 
			Because $\Abold$ is a symmetric indefinite matrix, any Krylov iterative method suitable for this class of matrices, such as
			MINRES \cite{Paige1975}, SYMMLQ \cite{Paige1975}, SQMR \cite{freund1994}, and GMRES, can be used.
			For successful convergence of the iterative method, one is often required to apply a good preconditioner.
      Even without any preconditioner, the saddle point system itself tends to be better conditioned for moderately large values of $\phi$ 
      than the Schur complement in the range-space method \cite{Thorne2009}. 
      If a good preconditioner is used, the full-space method is likely to have a better scalability than the range-space method,
      meaning that computational cost does not grow at an exponential rate as the problem size increases.
      Thus, many preconditioners have been developed recently. 
      Here we focus on a Schur complement-based preconditioner
      because the primary interest is to discuss the usage of a Schur complement factorization that will be
      introduced in Section~\ref{sec:exactrepresentation}.
			Murphy, et al.\ \cite{Murphy2000} have shown that if the following block diagonal preconditioner
			is used, then the preconditioned system has at most three nonzero distinct eigenvalues:
			\begin{equation}
			\label{eq:murphypreconditioner}
				\Pbold_{mgw} = \bmat{\Vbold    &  \zerobold   & \zerobold
							      		\\   \zerobold &  \phi \Vbold & \zerobold
									      \\   \zerobold &  \zerobold   & \Sbold    },
			\end{equation}
			where $\Sbold$ is the negative Schur complement defined in \eqref{eq:schurcomplementsystem}.
			This implies that the maximum number of iterations required for convergence 
			in a Krylov iterative method is three if $\Pbold_{mgw}$ is used as a preconditioner.
			However, applying this preconditioner has previously been considered impractical because it requires solving
			with $\Sbold$ in order to apply $\Pbold_{mgw}$. Even if it were practical to solve with $\Sbold$ efficiently
			and accurately, then the preferred approach would typically be to use the range-space method rather than to apply $\Pbold_{mgw}$ in the full-space method
      because the range-space method requires only one solve with $\Sbold$, while the full space method with $\Pbold_{mgw}$ as a preconditioner
      is likely to require more than one solve with $\Sbold$.
		  There has been some research done on finding a good approximation of the Schur complement 
      for use in a Schur complement-based preconditioner (e.g.,\ $\Pbold_{mgw}$) \cite{Rees2010, Pearson2012, Pearson2012siam}. 
      Particularly, Pearson and Wathen \cite{Pearson2012}	have developed the following approximation, which is regularization-robust:
			\begin{equation}
			\label{eq:pearsonapproximation}
			\begin{aligned}
				\Sbold_p &= (\Kbold+\frac{1}{\sqrt{\phi}}\Vbold)\Vbold^{-1}(\Kbold+\frac{1}{\sqrt{\phi}}\Vbold) \\
					  &= \Sbold + \frac{2}{\sqrt{\phi}}\Kbold,
			\end{aligned}
			\end{equation}
			and proved that the eigenvalues of $\Sbold_p^{-1}\Sbold$ are between $\frac{1}{2}$ and $1$ regardless of the mesh size $h$ 
			and $\phi$. This implies theoretically that a Krylov iterative method must
      converge in $O(1)$ iterations regardless of $h$ and $\phi$.
			They have used an algebraic multi-grid method to solve with the factor $\left (\Kbold+\frac{1}{\sqrt{\phi}}\Vbold \right )$ in their numerical examples
      and demonstrated that the number of iterations required for convergence is indeed $O(1)$.
      However, in the problems they solve, the condition number of the Schur complement is small enough
      so that the approximation to the Schur complement behaves well.
			The factorization introduced in Section~\ref{sec:exactrepresentation} has a similar form to \eqref{eq:pearsonapproximation}.
			However, it is an exact representation of $\Sbold$ and therefore enables the application of the range-space method.

	\section{A ``practical" factorization of the Schur complement}
	\label{sec:exactrepresentation}
		The Schur complement $\Sbold$ in \eqref{eq:schurcomplementsystem} can be factored into the following form:
		\begin{equation}
		\label{eq:factorizationofschur}
			\Sbold_c = (\Kbold+\frac{1}{i\sqrt{\phi}}\Vbold)\Vbold^{-1}(\Kbold-\frac{1}{i\sqrt{\phi}}\Vbold),
		\end{equation}
    where $i=\sqrt{-1}$.
		The form of this factorization is similar to $\Sbold_p$ in \eqref{eq:pearsonapproximation} in the sense that
		the first and third factors are some linear combinations of $\Kbold$ and $\Vbold$, and the middle factor is $\Vbold^{-1}$.
    This suggests that the same methods for solving with the factors in $\Sbold_p$ may also be used to solve with factors
    in $\Sbold_c$. However, there are some important differences between $\Sbold_p$ and $\Sbold_c$. The first and third factors of $\Sbold_p$ are the same, but 
		the corresponding factors of $\Sbold_c$ are not. Consequently, if direct solvers were to be used for each factor,
    then only one factorization would be required for $\Sbold_p$, while two factorizations of two different systems
    would be required for $\Sbold_c$. Additionally, $\Sbold_c$ introduces complex numbers, but all the elements in $\Sbold_p$
		are real. Thus, one solve with $\Sbold_c$ requires four times more storage and floating point operations than one solve with $\Sbold_p$.
    In spite of the disadvantages of applying $\Sbold_c$, it is an exact representation of $\Sbold$, and this permits 
    use of the range-space method where only one solve with $\Sbold_c$ is sufficient
    to obtain a solution to \eqref{eq:discretizedoptimalheatcontrol}.
    On the other hand, many solves with $\Sbold_p$ are required because $\Sbold_p$ is an approximation to $\Sbold$, so 
    it can only be used as a preconditioner. This discussion leads to the following two extreme cases
    when dealing with the Schur complement:
    \begin{itemize}
      \item If a direct method is the preferred choice for solving the Schur complement 
            (i.e., a factorization of a given system is required), 
            it would be advantageous to apply the direct method to $\Sbold_p$ 
            assuming that the dominating computational cost occurs in factorization of the system.
      \item If an iterative method is the only option and the range-space method is applicable, then
            it would be preferable to apply the iterative method to $\Sbold_c$.      
    \end{itemize}
    Between these two extreme cases, a choice must be made depending on the characteristics of
    the problem and a numerical solver. For example, if the domain of a problem is complex, then the computational overhead incurred by 
    applying $\Sbold_c$ rather than $\Sbold_p$ will be substantially diminished since in this case the factors of both $\Sbold_c$ and $\Sbold_p$ will be complex. 
    Such problems include any frequency domain analyses with damping in structural or 
    acoustic problems. On the other hand, if a problem requires
    many solves with a Schur complement and with multiple right-hand sides
    and is small enough to permit a direct solver to be used for the factors of $\Sbold$,
    applying $\Sbold_p$ is favorable because
    direct methods are more efficient than iterative methods in general for multiple right-hand sides.
    If a domain decomposition method---in which both direct and iterative methods are used---is chosen,
    then one needs to examine the costs of each component of the computation (e.g., data structure,
    building and storing the operator, and solving with the operator).
    In Section~\ref{sec:numericalresults}, numerical results for solving with $\Sbold_c$
    are shown using the domain decomposition solver FETI-DPH \cite{Farhat2005}.
			
	\section{FETI}
	\label{sec:FETI}
		In order to be self-contained, the FETI method \cite{farhat1991} and two of its variants 
		(FETI-DP \cite{Farhat2001} and FETI-DPH \cite{Farhat2005}) are explained in this section. 
		For a more detailed description, see \cite{farhat1991, Farhat2001, Farhat2005, Avery2009}.
		The FETI method was developed in order to solve, in parallel, the following linear system of equations
		arising from the finite element discretization of a linear elasticity PDE:
		\begin{equation}
		\label{eq:linearsystem}	
			\Kbold\ybold = \fbold,
		\end{equation}
		where $\Kbold$ is the stiffness matrix (or a linear combination of the stiffness matrix and mass matrix), 
    $\ybold$ is the vector of unknown displacements, and $\fbold$ is an external force term.
		Solving \eqref{eq:linearsystem} is equivalent to minimizing a quadratic function:
    \begin{equation}
    \label{eq:feti_qp}
		\begin{aligned}
			\underset{\ybold}{\text{minimize}} & & \frac{1}{2} \ybold^T \Kbold \ybold - \ybold^T \fbold.
		\end{aligned}
    \end{equation}
		\begin{figure}[t]
    \centering
			\includegraphics[scale=0.4]{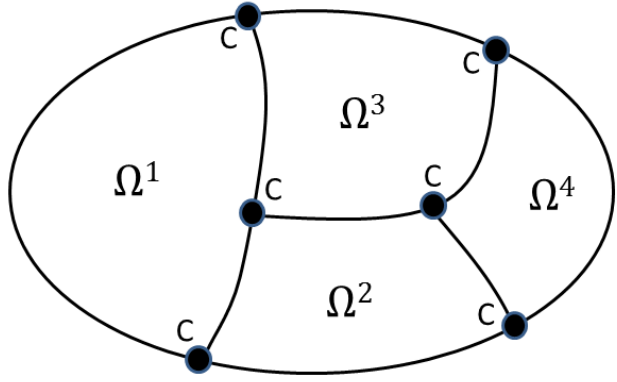}
		\caption{Illustration of 4 non-overlapping subdomains. The sample corner nodes are denoted by C. 
             FETI-DP enforces continuity on these corner nodes.}
		\label{fig:subdomains}
		\end{figure}
		The objective function in \eqref{eq:feti_qp} often represents a physical quantity (e.g., energy in linear elasticity).
		In domain decomposition methods, the spatial domain $\Omega$ is divided into $N_s$ either overlapping or non-overlapping subdomains.
		The FETI method divides $\Omega$ into $N_s$ non-overlapping subdomains as illustrated in Figure \ref{fig:subdomains} 
    and minimizes a local function in each subdomain
		(i.e., minimize ${\ybold^s}^T\Kbold^s\ybold^s - {\ybold^s}^T \fbold^s$ in $\Omega^s$ 
    where the superscript ``s'' designate the restrictions to the specific subdomain).
    In order to be equivalent to \eqref{eq:feti_qp},
		the following additional continuity condition on $\ybold^s$ between interfaces must be imposed:
		\begin{equation}
    \label{eq:continuitycondition} 
      \sum_{s=1}^{N_s} \Bbold^s\ybold^s = \zerobold,
    \end{equation}
		where $\Bbold^s$ is a signed Boolean matrix. 
		Hence one can write the following constrained quadratic programming equivalent to \eqref{eq:feti_qp}:
		\begin{equation}
    \label{eq:feti_cqp}
		\begin{array}{cl}
      \underset{\ybold^s; s=1,\ldots,N_s}{\text{minimize}} & \sum_{s=1}^{N_s} 
              \left (\frac{1}{2} {\ybold^s}^T \Kbold^s {\ybold^s} - {\ybold^s}^T \fbold^s \right )
			\\ \text{subject to} & \sum_{s=1}^{N_s} \Bbold^s\ybold^s = \zerobold. 
		\end{array}
		\end{equation}
		The Lagrangian for \eqref{eq:feti_cqp} is 
		\begin{equation}
		\label{eq:lagrangian_cqp}
			L(\ybold^s, \lambdabold; s=1,\ldots,N_s) = \sum_{s=1}^{N_s} 
             \left (\frac{1}{2}{\ybold^s}^T \Kbold^s {\ybold^s} - {\ybold^s}^T \fbold^s \right ) + \lambdabold^T \sum_{s=1}^{N_s} \Bbold^s\ybold^s.
		\end{equation}
		For simplicity of notation, we will omit the writing of $s=1,\ldots,N_s$. The variable $\ybold^s$ 
		either means an individual subvector on $\Omega^s$ or a set of subvectors over all the subdomains.
		Because \eqref{eq:feti_cqp} is convex and only linear equality constraints are present,
		Slater's condition holds. Therefore strong duality also holds.
		Thus, an optimal solution ${\ybold^s}^*$ to \eqref{eq:feti_cqp} can be obtained 
		from a saddle point (${\ybold^s}^*$,$\lambdabold^*$) to \eqref{eq:lagrangian_cqp}, that is,
		\begin{equation}
		\label{eq:saddlepoint}
			\sup_{\lambdabold} \inf_{\ybold^s} L(\ybold^s,\lambdabold) = L({\ybold^s}^*,\lambdabold^*) = 
            \inf_{\ybold^s} \sup_{\lambdabold} L(\ybold^s, \lambdabold).
		\end{equation}
		Statement \eqref{eq:saddlepoint} says that one can find an optimal solution ${\ybold^s}^*$ in two ways.
    The first way is to minimize $L(\ybold^s,\lambdabold)$ with respect to $\ybold^s$, then maximize the resultant with respect to $\lambdabold$.
    The second way is to maximize $L(\ybold^s,\lambdabold)$ with respect to $\lambdabold$, then minimize the resultant with respect to $\ybold^s$.
		The FETI method takes the first approach where the Lagrange dual function $g(\lambdabold) = \inf_{\ybold^s} L(\ybold^s, \lambdabold)$
		is obtained and then maximized in order to obtain the Lagrange multipliers first.
    This approach is attractive if the number of constraints (i.e., interface continuity condition) is small, 
    which is expected to be the case in general.

		For a fixed $\lambdabold$, $L(\ybold^s, \lambdabold)$ is separable in $\ybold^s$ 
		(i.e., $\left (\frac{1}{2}{\ybold^s}^T \Kbold^s {\ybold^s} - {\ybold^s}^T \fbold^s + \lambdabold^T \Bbold^s\ybold^s \right )$ on $\Omega^s$). 
		Thus, the Lagrange dual function can be obtained by	solving $\Kbold^s\ybold^s = \fbold^s - \Bbold^s\lambdabold$ on each subdomain $\Omega^s$.
    In the case of $\Kbold$ being a stiffness matrix as in \eqref{eq:linearsystem},
		$\Kbold^s$ is symmetric positive definite or semidefinite. If $\Kbold^s$ is positive semidefinite, then
		it is necessary to explicitly ensure that $\Kbold^s\ybold^s = \fbold^s - \Bbold^s\lambdabold$ is compatible 
    (i.e., ${\Rbold^s}^T(\fbold^s-{\Bbold^s}^T\lambdabold) = \zerobold$
		where the columns of $\Rbold^s$ span the left null space of $\Kbold^s$). Otherwise,
		the minimum of $\left (\frac{1}{2}{\ybold^s}^T \Kbold^s {\ybold^s} - {\ybold^s}^T \fbold^s + \lambdabold^T \Bbold^s\ybold^s \right )$ in terms of $\ybold^s$ becomes $-\infty$
		and the objective value of primal problem \eqref{eq:feti_cqp} is also $-\infty$ by strong duality.
		However, this is not a physical solution. Thus, we restrict ourselves to the case when $g(\lambdabold) > -\infty$.
		Finally, the Lagrange dual function $g(\lambdabold) = \inf_{\ybold^s} L(\ybold^s, \lambdabold)$ is defined as
		\begin{equation}
			g(\lambdabold) = \left \{ \begin{array}{ll} -\frac{1}{2}	\lambdabold^T \Fbold \lambdabold + \dbold^T\lambdabold - \cbold	
                                 & \text{if } {\Rbold^s}^T(\fbold^s-{\Bbold^s}^T\lambdabold) = \zerobold \text{ for } \forall s, 
                              \\ -\infty              & \text{otherwise}, \end{array} \right .
		\end{equation}
		where 
		\begin{equation}
		\label{eq:defineFdc}
		\begin{aligned}			
			\Fbold = \sum_{s=1}^{N_s} \Bbold^s{\Kbold^s}^+{\Bbold^s}^T, \quad
			\dbold = \sum_{s=1}^{N_s} \Bbold^s{\Kbold^s}^+\fbold^s, \quad
			\cbold = \frac{1}{2} \sum_{s=1}^{N_s} {\fbold^s}^T{\Kbold^s}^+\fbold^s,
		\end{aligned}
		\end{equation}
    and ${\Kbold^s}^+$ is a pseudo-inverse of $\Kbold^s$. This defines the following dual problem to \eqref{eq:feti_cqp}:
		\begin{equation}
    \label{eq:dual_cqp}
		\begin{array}{cl}
      \underset{\lambdabold}{\text{maximize}} & -\frac{1}{2}  \lambdabold^T \Fbold \lambdabold + \dbold^T\lambdabold 
      \\ \text{subject to} & {\Rbold^s}^T(\fbold^s-{\Bbold^s}^T\lambdabold) = \zerobold \quad \text{for } \forall s. 
		\end{array}
    \end{equation}
		The FETI method solves \eqref{eq:dual_cqp} with a Preconditioned Conjugate Projected Gradient (PCPG) algorithm.
		The FETI formulation above is suitable for parallel processing. Any subdomain level computations (e.g., 
		factorization or computation with $\Kbold^s$ and $\fbold^s$) can be assigned to an individual process. 
		The size of $\lambdabold$ is the total number of degrees of freedom restricted to the interfaces between subdomains.
		Within PCPG, $\lambdabold$ is projected onto the domain of feasibility, 
		which in turn reduces the effective dimension of the problem.  
		Indeed the FETI method equipped with the Dirichlet preconditioner
    is proven numerically scalable with respect to both problem size and number of subdomains.
    The Dirichlet preconditioner is defined as 
    \begin{equation}
    \label{eq:dirichletpreconditioner}
      \Pbold^{-1} = \Wbold\sum_{s=1}^{N_s} \Bbold^s\bmat{\zerobold & \zerobold \\ \zerobold & \Sbold^s_{bb}}{\Bbold^s}^T\Wbold,
    \end{equation}
    where 
    \begin{equation}
    \label{eq:WS}
       \Wbold = \left ( \sum_{s=1}^{N_s} \Bbold^s {\Bbold^s}^T \right )^+, \quad \Sbold^s_{bb} = \Kbold_{bb}^s - {\Kbold_{ib}^s}^T{\Kbold_{ii}^s}^{-1}\Kbold_{ib}^s,
    \end{equation}
    and the subscript $i$ and $b$ denote 
    subdomain internal degrees of freedom and interface degrees of freedom, respectively.
		If it is applied in the conjugate gradient algorithm with the projected gradient \cite{gill1974numerical}
    to second-order elliptic problems,
		the condition number $\kappa$ of the interface problem \eqref{eq:dual_cqp} is approximately \cite{Mandel1996}
		\begin{equation}
		\label{eq:conditionnumberbound}
			\kappa = O(1+\log^m(H/h)), \quad m \leq 3.
		\end{equation}
		However, the first-generation FETI method is not numerically scalable for fourth-order plate and shell problems. 
		This leads to FETI-DP. FETI-DP is one variant considered to be the third-generation FETI method.
    FETI-DP enforces continuity at some interface corner nodes at each iteration (see Figure \ref{fig:subdomains}).
    An extra coarse problem needs to be solved \cite{Farhat2001} and
    in each subdomain the remaining subdomain stiffness matrix $\Kbold^s_{rr}$ 
    (i.e., that excluding the degrees of freedom in the corner nodes at which the continuity is enforced)
    becomes positive definite. Consequently, the corresponding dual interface problem becomes an unconstrained QP. 
    This treatment helps to achieve numerical scalability for fourth-order elliptic problems. 
    The continuity constraints can be augmented by additional constraints that are enforced exactly 
    throughout the iterations in FETI-DP in order to accelerate the convergence.
    This augmentation procedure results in the augmented coarse problem in FETI-DP \cite{lesoinne200319}.
    A standard augmentation procedure uses the edge-based rigid body modes (rotational and/or translational) \cite{Farhat2005}.
    The positive definiteness of $\Kbold^s_{rr}$ is not guaranteed when $\Kbold$ in \eqref{eq:linearsystem} is different
    from the stiffness matrix. For example, for a Helmholtz problem or frequency response elastodynamic problem, 
    $\Kbold$ in \eqref{eq:linearsystem} becomes
    \begin{equation}
    \label{eq:helmholtzoperator}    
      \Zbold \equiv \Kbold - k^2\Mbold + i\Cbold,
    \end{equation} 
    where $\Cbold$ is a symmetric matrix that arises from the discretization of an absorbing boundary condition, 
    and $k>0$ is a frequency (or a wave number for acoustic scattering problems). 
    In this case, depending on the value of $k$, $\Zbold$ may become indefinite.
    The same difficulty is encountered when solving with complex factors in \eqref{eq:factorizationofschur}.
    In order to resolve this difficulty, two special treatments were required.
    \begin{itemize}
      \item The first treatment is required to deal with indefinite matrices.
            A solver that is suitable for indefinite matrices must be applied.
            Such an iterative solver includes GMRES \cite{Saad1986} for a general square matrix
            and MINRES \cite{Paige1975} for a symmetric indefinite matrix.
            FETI-DPH (`H' stands for Helmholtz) \cite{Farhat2005} uses GMRES to deal with indefiniteness.
      \item The second treatment is required for improved convergence when the wave number of frequency is large.
            A particular augmentation used in FETI-DPH that accelerates convergence for the Helmholtz or frequency response elastodynamic problem 
            is the free-space solutions of the corresponding equations, which are plane waves.
    \end{itemize}
    Note that the complex factors in Section \ref{eq:factorizationofschur} are similar to \eqref{eq:helmholtzoperator}, 
    but not identical, and not surprisingly the plane wave augmentation was found to not provide the same beneficial effects as in
    the Helmholtz problem. Thus, the edge-based rigid body modes
    (both rotational and translational) augmentation is applied in the next section for two numerical examples.
    Strictly speaking, although FETI-DPH normally refers to the variant of FETI-DP that uses both GMRES and plane wave augmentation, 
    in what follows, FETI-DPH refers to a FETI-DP solver incorporating GMRES and edge-based rigid body modes.


	\section{Numerical Results}
	\label{sec:numericalresults}
    Two pedagogical examples are considered: linear static heat control of a 2D square plate and 
    structure control of a 3D solid cantilever. The heat control problem in Section \ref{sec:thermalproblem}
    is taken identical to the thermal problem solved in the paper by Rees et al. \cite{Rees2010}.
    The Schur complement that arises in this thermal problem is sufficiently well-conditioned for the range of
    $\phi$ values and mesh size $h$ considered here for the range-space method to be applicable and therefore it is the only method used.
    The cantilever control problem in Section \ref{sec:cantileverproblem} is very flexible with material properties of rubber.
    For a relatively large $\phi$ value, the range-space method fails to converge to an accurate solution.
    Thus, the range-space method is used for small values of $\phi$ only, and
    the full-space method is used for large values of $\phi$. In the case of the full-space method, 
    both the approximate Schur complement proposed by Pearson and Wathen \cite{Pearson2012} and
    the complex factorization introduced in this paper are used in a Schur-complement based preconditioner 
    and compared in terms of computational time and iteration counts.
  
    It is well established that augmentation is a beneficial feature that improves the performance of FETI. 
    However, the optimality of the augmentation's performance depends on the nature of the problem. 
    Because neither the FETI solver nor any of its variants has ever been applied to a system of the form
    $\Kbold\pm\frac{1}{i\sqrt{\phi}}\Vbold$,
    the numerical performance of FETI-DPH with and without augmentation is compared in this section.
    All the simulations were run on a heterogeneous Linux cluster with 2.6 GHz hexa-core Westmere processors 
    (32 blades, 2 processors per blade, QDR Infiniband interconnect, 24GB/blade) 
    and 2.6 GHz octa-core Sandybridge processors (4 blades, 2 processors per blade, FDR Infiniband interconnect, 256 GB/node).

  \subsection{Thermal problem}
  \label{sec:thermalproblem}
\begin{figure}[t]
    \centering
			\includegraphics[scale=0.25]{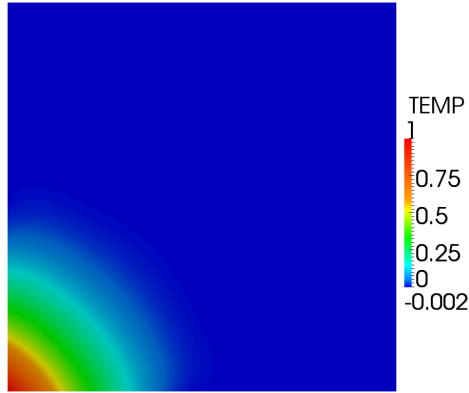}
		\caption{Target temperature defined as in \eqref{eq:ex1target}. }
		\label{fig:targettemperature}
		\end{figure}

\begin{figure}[t]
		\begin{center}
    $\begin{array}{c@{\hspace{0in}}c}
			\includegraphics[scale=0.18]{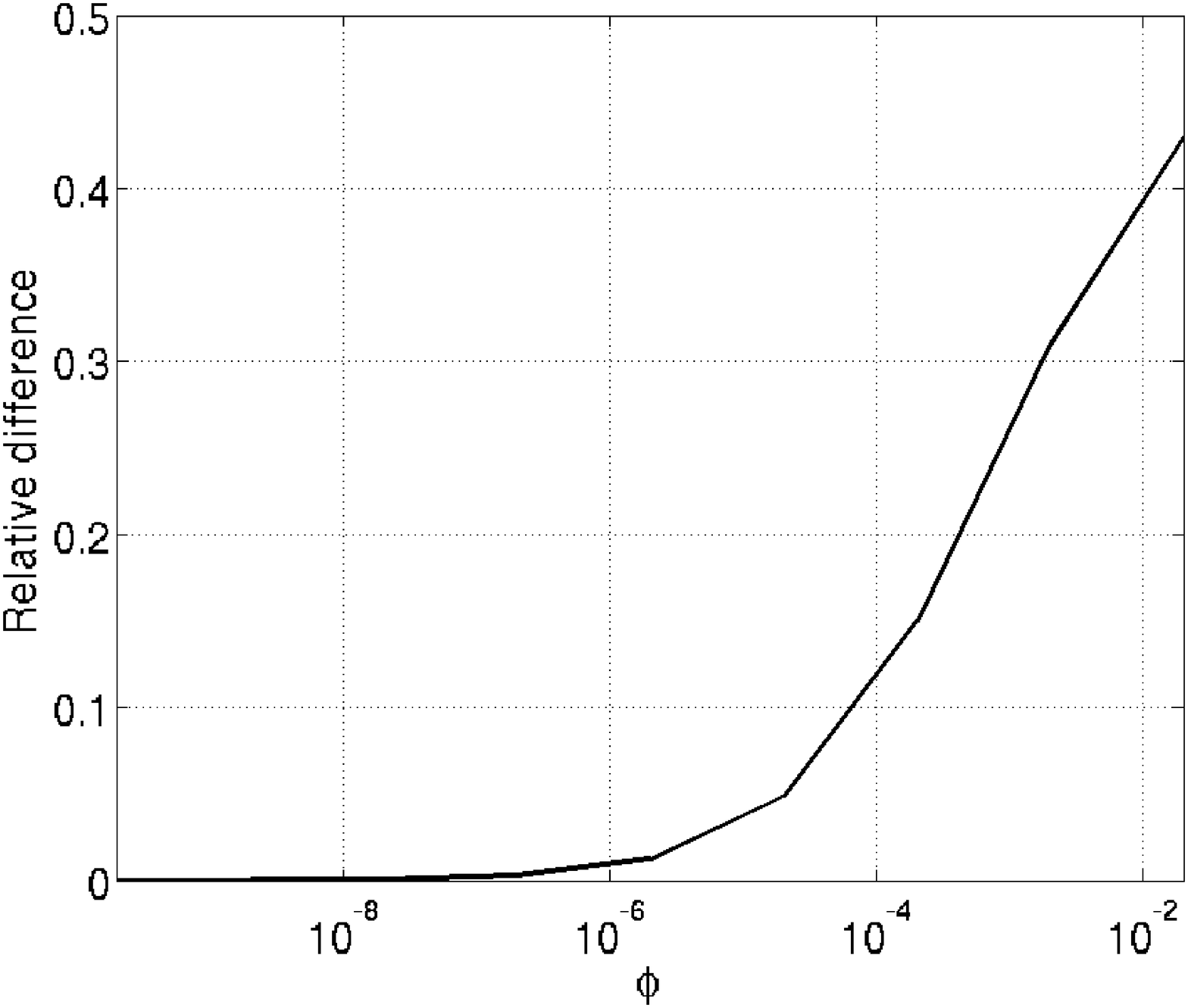} &
      \includegraphics[scale=0.18]{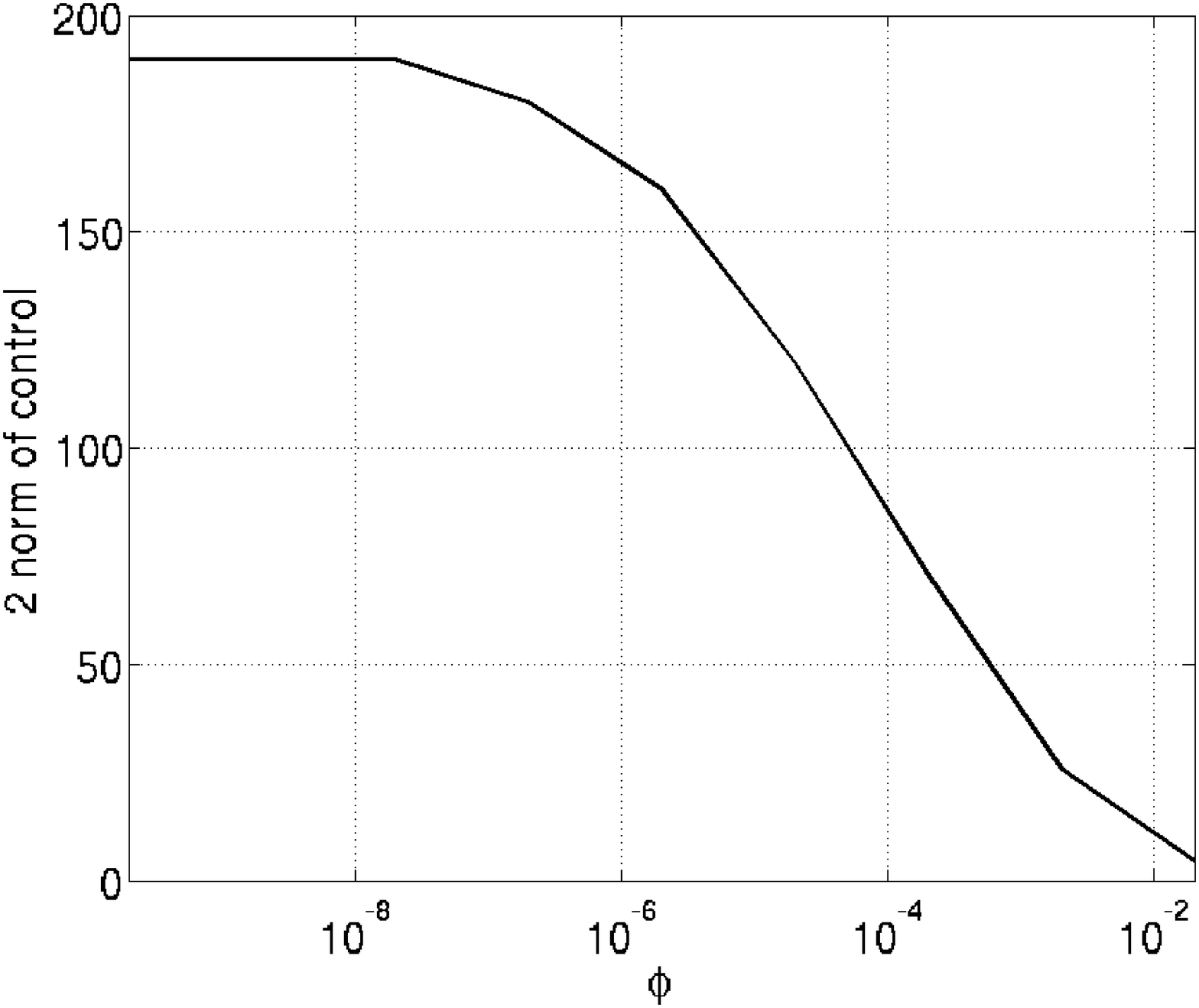} \\
    \mbox{\bf (a) } & \mbox{\bf (b)}
    \end{array}$
    \end{center}
		\caption{\textbf{(a)} Graph of $\phi$ vs relative difference between temperature solution and target temperature.
             \textbf{(b)} Graph of $\phi$ vs norm of control solution.}
		\label{fig:phivsnormdifftarget}
		\end{figure}

		A linear static heat control problem is solved. 
		This example is identical to example $5.1$ in \cite{Rees2010} by Rees et al.
		The same objective function and PDE constraints are used as in \eqref{eq:continuousoptimalheatcontrol}, whose continuous
		formulation is reproduced here for better access:
    \begin{equation}
    \label{eq:heatconduction}
    \begin{aligned}
			\underset{y,u}{\text{minimize}} & & F(y,u) &:= \frac{1}{2} \int_\Omega(y - \bar{y})^2 dx + \frac{\phi}{2} \int_\Omega u^2 dx
			\\ \text{subject to} & & -\nabla^2y &= u \text{ on } \Omega 
			\\ & & y &= y_c   \text{ on } \Gamma_g.
    \end{aligned}
    \end{equation}
    The domain $\Omega$ is $[0,1]^2 \subset \mathbb{R}^2$, which is a unit square plate,
    whose heat conductivity is unity. The target temperature $\bar{y}$ is defined as
    \begin{equation}
    \label{eq:ex1target}    
      \bar{y} = \left\{ \begin{array}{ll} (2x_1-1)^2(2x_2-1)^2 & \quad \text{if } (x_1,x_2) \in [0,\frac{1}{2}]^2, \\
                                          0                    & \quad \text{otherwise,} \end{array} \right.
    \end{equation}
    which is illustrated in Figure~\ref{fig:targettemperature}. 
		The boundary condition $y_c$ is defined as 
    \begin{equation}
    \label{eq:ex1boundarycondition}   
      y_c = \bar{y} \quad \text{ on } \Gamma_g = \partial \Omega = \{(x_1,x_2) \mid x_1 \in\{0,1\}, x_2 \in\{0,1\} \}.
    \end{equation}
    The optimal control problem \eqref{eq:heatconduction} tries to find a temperature $y$ that is close to
    the target temperature $\bar{y}$ by controlling heat $u$. How close $y$ can be to $\bar{y}$ is
    determined by the parameter $\phi$. As $\phi$ decreases, $y$ is expected to approach $\bar{y}$,
    but $\|u\|$ is also expected to increase. This makes sense because the objective function
    value is not sensitive to the second term if the parameter $\phi$ is small and the first
    term dominates the objective function value.
		The results shown in Figure~\ref{fig:phivsnormdifftarget} confirm the validity of these expectations. 
    For this particular example, according to Figure~\ref{fig:phivsnormdifftarget}(a),
    one needs $\phi$ less than $2\times10^{-6}$ in order to match the target temperature to within $1\%$.
    According to Figure~\ref{fig:phivsnormdifftarget}(b), however,
		one needs to set $\phi$ greater than $2\times10^{-4}$ in order for the norm of control solution not to exceed $100$. 

    \begin{figure}[th]
    \begin{center}
    $\begin{array}{c@{\hspace{0.1in}}c@{\hspace{0.1in}}c}
      \includegraphics[scale=0.13]{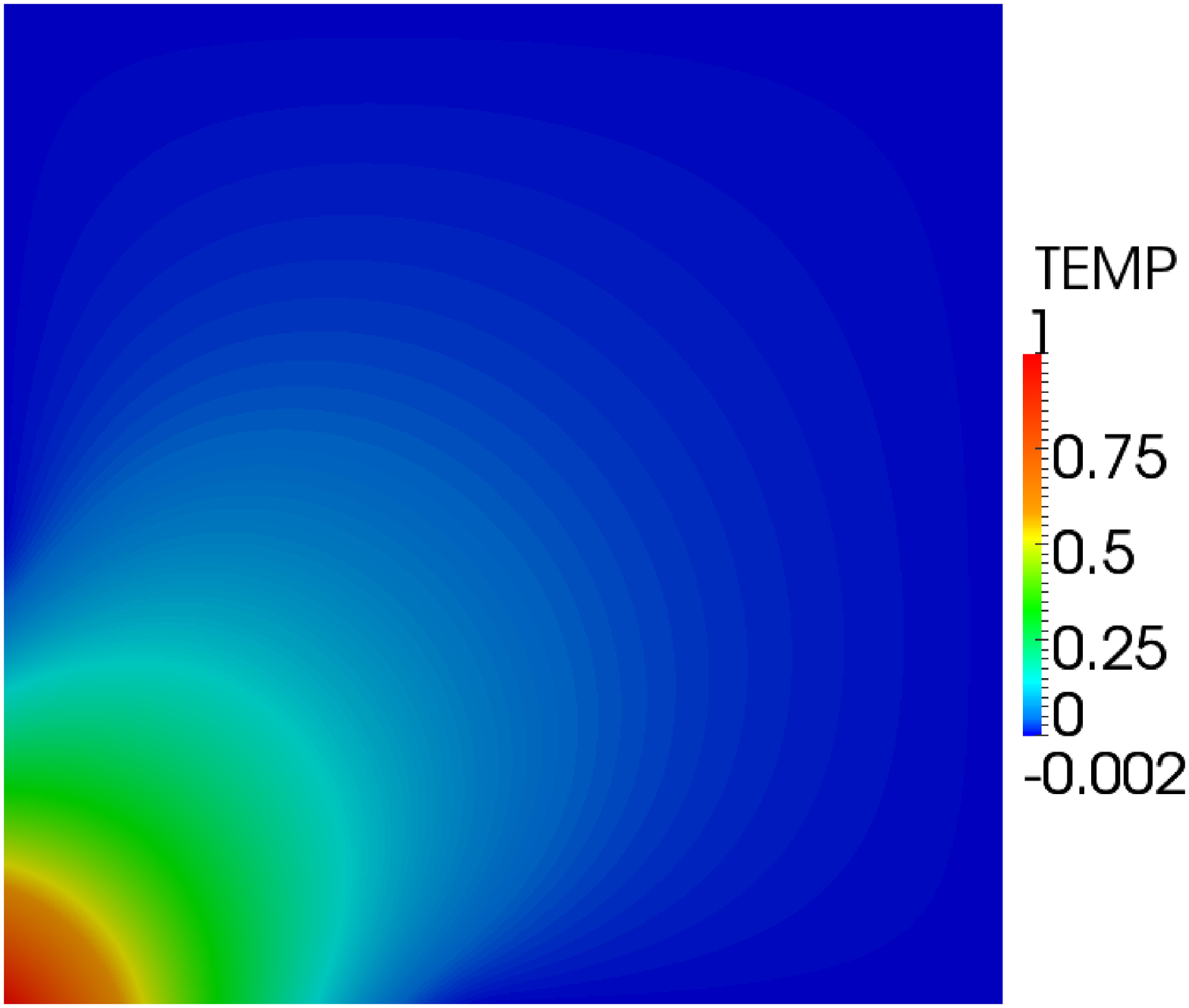} &
      \includegraphics[scale=0.13]{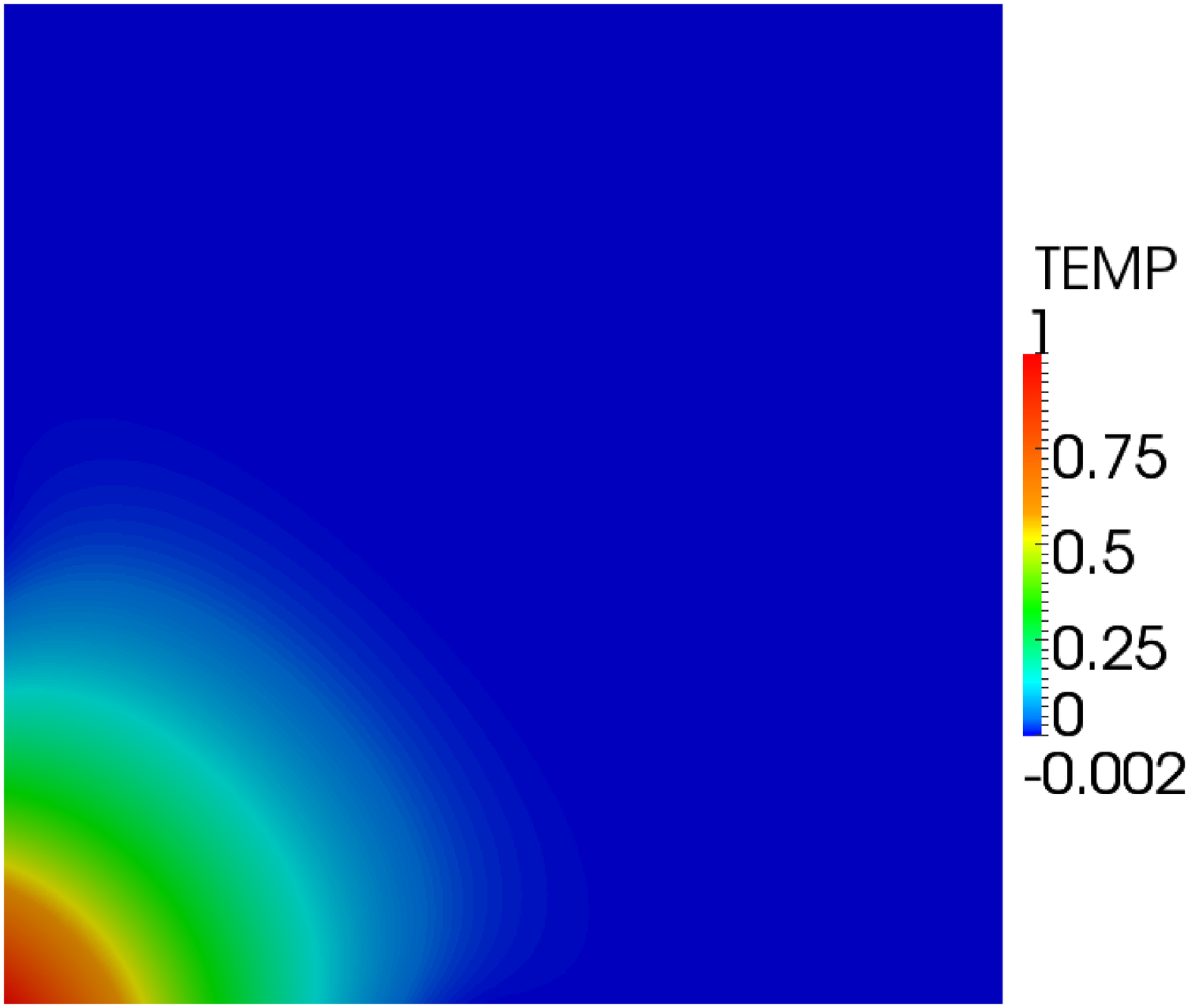} &
      \includegraphics[scale=0.13]{temp5.eps} \\ 
    	\includegraphics[scale=0.13]{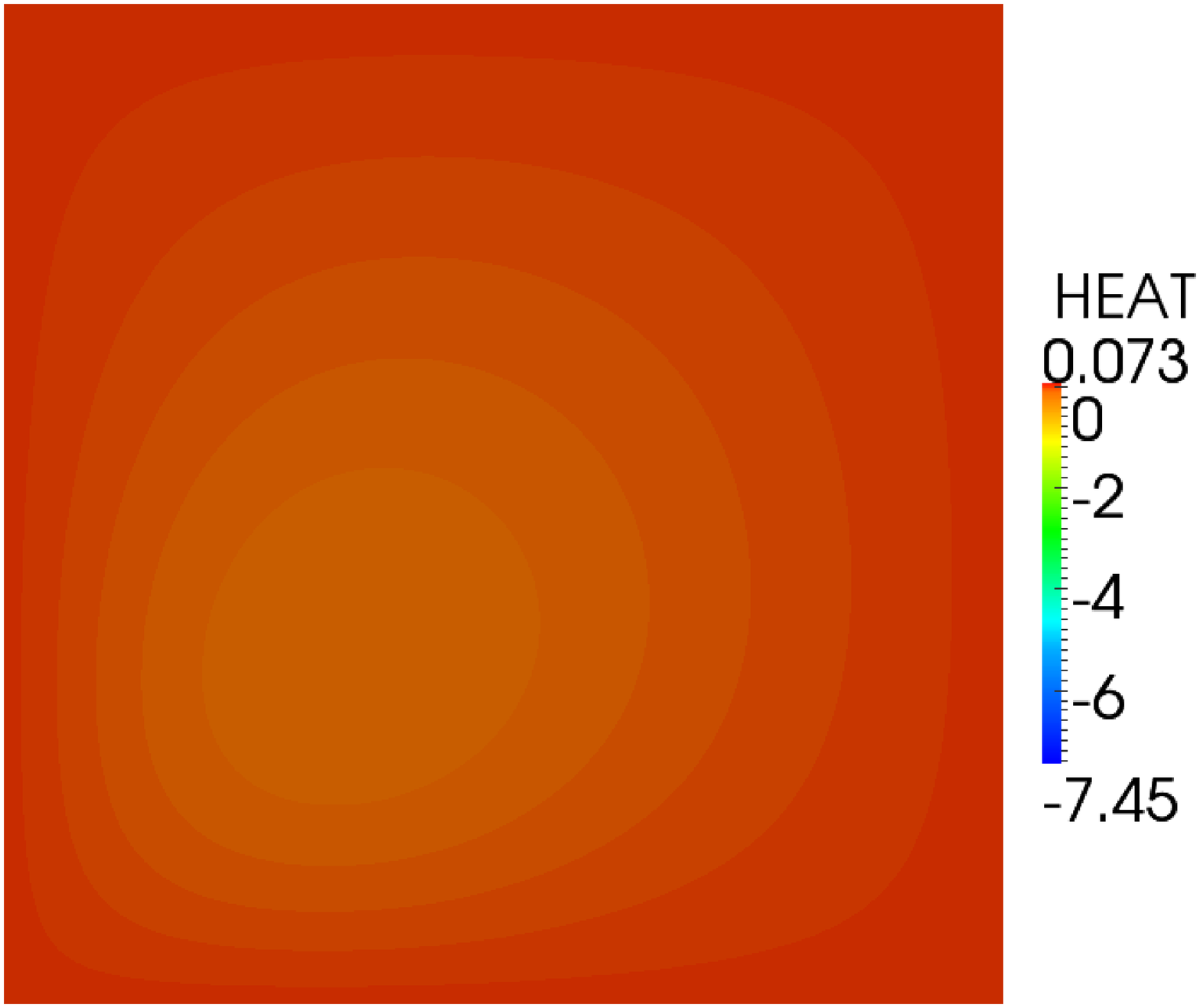} &
    	\includegraphics[scale=0.13]{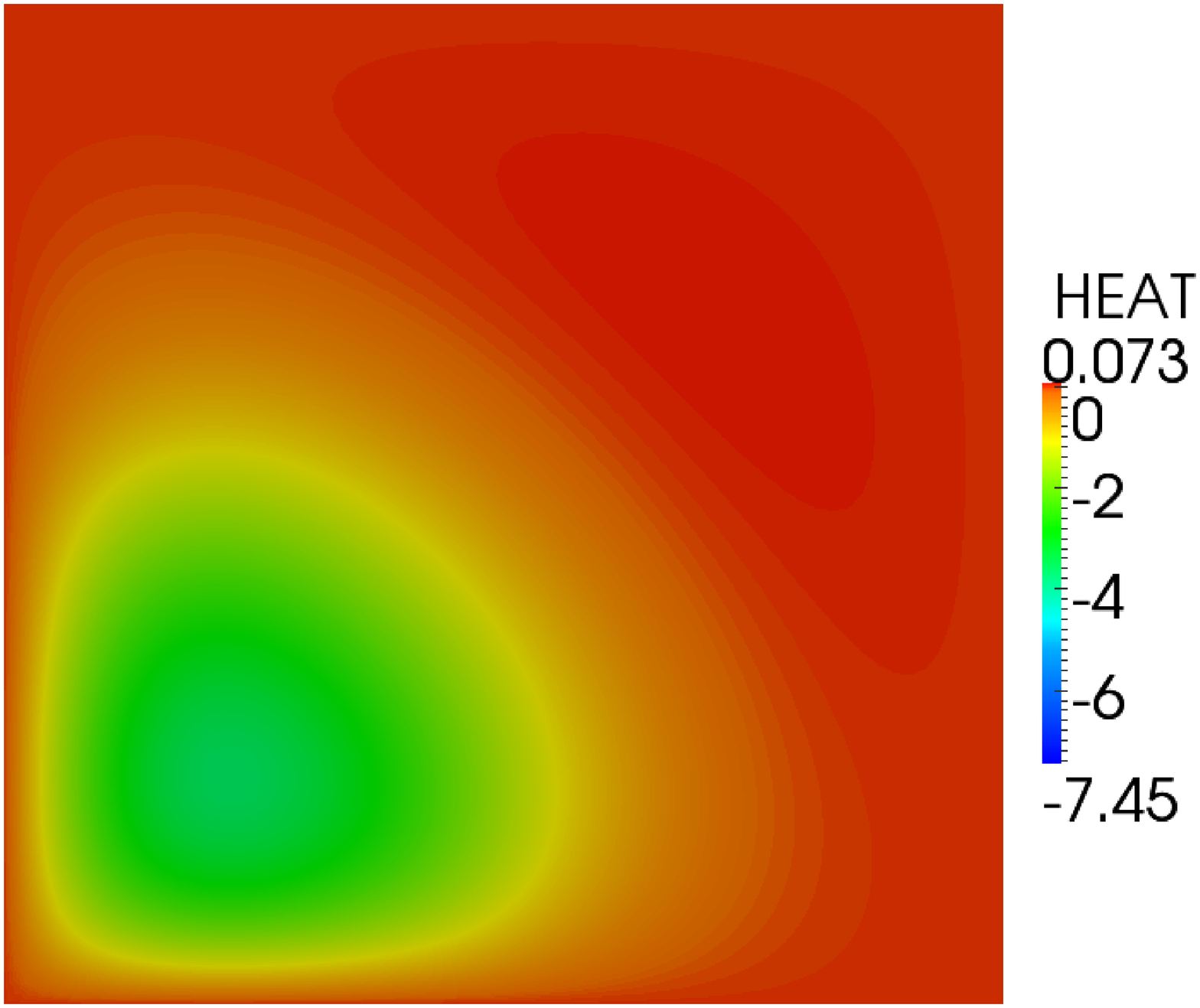} &
      \includegraphics[scale=0.13]{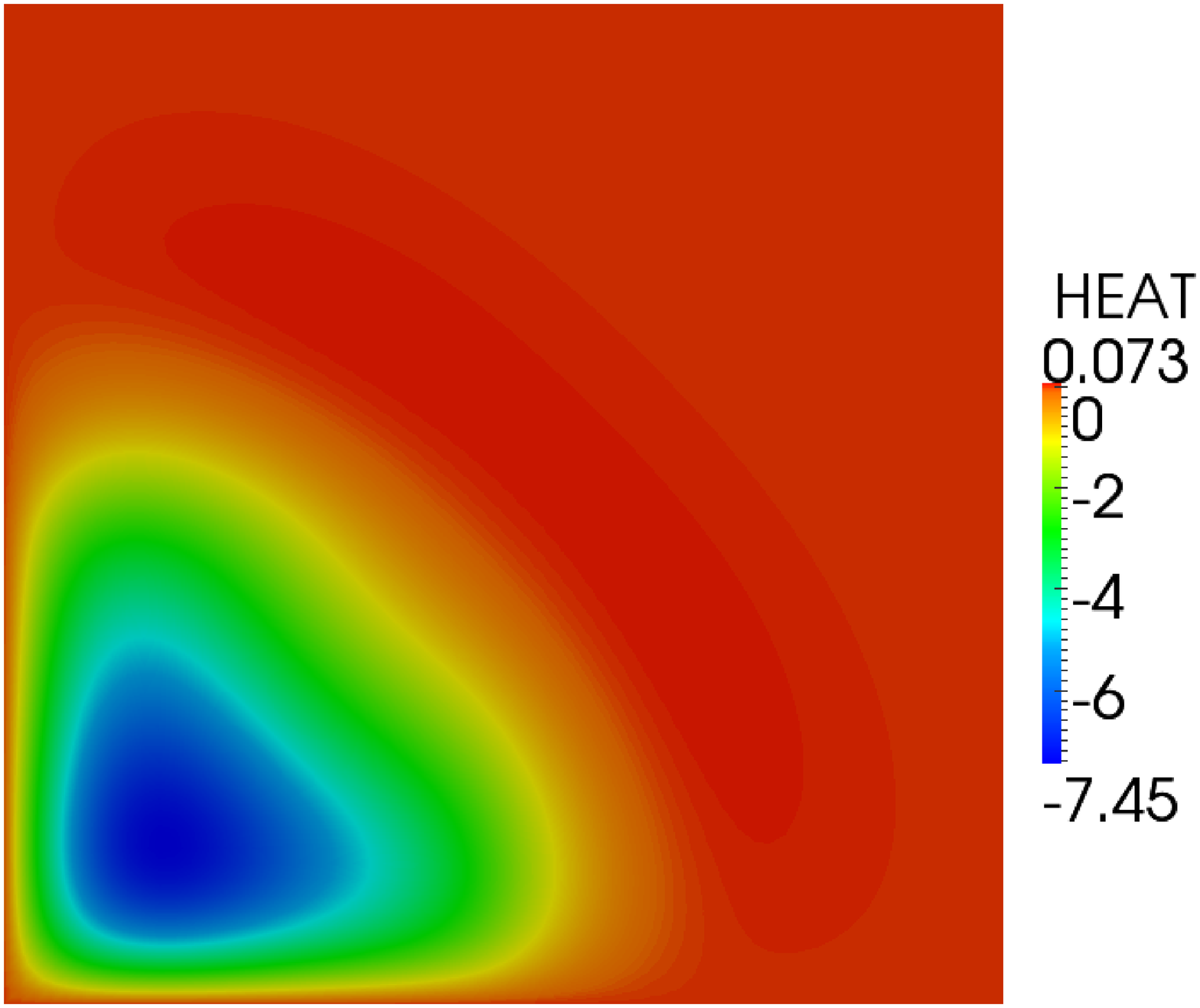} \\ 
    \mbox{\bf (a) $\phi = 2\times10^{-2}$} & \mbox{\bf (b) $\phi = 2\times10^{-4}$} & \mbox{\bf (c) $\phi = 2\times10^{-5}$}
    \end{array}$
    \end{center}
    \caption{Temperature (the upper three figures) and 
             heat distribution solutions (the lower three figures) for various $\phi$ values.} 
    \label{fig:ex1tempheatforvariousregularization}
    \end{figure}

    The same effects can be demonstrated visually in Figure~\ref{fig:ex1tempheatforvariousregularization}.
    Temperature and heat distributions for various $\phi$ values
    are shown in Figures~\ref{fig:ex1tempheatforvariousregularization}(a)-(c).
    Note that for $\phi = 0.02$, the heat is almost zero everywhere and 
    the corresponding temperature distribution (induced mainly by the boundary condition $y_c$) 
    is slightly different from the target temperature (i.e., Figure~\ref{fig:targettemperature}).
    The target temperature can be matched more precisely if a smaller $\phi$ is used.
    Figures~\ref{fig:ex1tempheatforvariousregularization}(b) and (c) show temperature distributions
    that are closer to the target temperature. They are produced by using smaller $\phi$ values
    (i.e., $\phi = 2\times10^{-4}$ and $2\times10^{-5}$). Note that the corresponding heat distributions
    show noticeable non-zero heat at the left bottom of the domain. 
    As $\phi$ decreases, a sharper heat gradient is visible near the boundary.

    \begin{table}[th]
    \caption{FETI-DPH's $\phi$ dependency for the heat control problem.}
    \centering
    \footnotesize
    \begin{tabular}{|c||c|c|c||c|c|c|}
         \hline   $\phi$ &  -  &  +  & CPU  &  RAUG- & RAUG+ & CPU
      \\ \hline  
         \hline   2E-2   & 33  & 32  & 9.4  &  13    & 13    & 6.2 
      \\ \hline   2E-3   & 36  & 35  & 10.2 &  13    & 13    & 6.3
      \\ \hline   2E-4   & 39  & 38  & 10.4 &  13    & 13    & 6.2
      \\ \hline   2E-5   & 38  & 37  & 10.4 &  13    & 13    & 6.3
      \\ \hline   2E-6   & 34  & 34  & 10.5 &  13    & 13    & 6.3
      \\ \hline   2E-7   & 34  & 34  & 9.7  &  13    & 13    & 6.2
      \\ \hline   2E-8   & 38  & 38  & 10.6 &  13    & 13    & 6.1
      \\ \hline   2E-9   & 42  & 41  & 11.2 &  14    & 14    & 6.5
      \\ \hline   2E-10  & 47  & 44  & 12.1 &  18    & 17    & 7.0
      \\ \hline   2E-11  & 56  & 52  & 14.3 &  22    & 21    & 7.7
      \\ \hline   2E-12  & 69  & 62  & 15.8 &  28    & 26    & 8.8
      \\ \hline
    \end{tabular}
    \label{ta:FETIDPphidependence}
    \end{table}
		
		Table~\ref{ta:FETIDPphidependence} shows FETI-DPH's dependence on $\phi$ and
    the effects of augmentation in the FETI-DPH solver.
		The regularization parameter $\phi$ varies from $2\times10^{-12}$ to $2\times10^{-2}$.
    The mesh size is $2^{-11}$ and the convergence threshold is $10^{-10}$.
    The size of each subdomain is $2^{-5}$ and the number of processes is $48$.
    The second and third columns show the number of FETI-DPH iterations without any augmentation, and
    the fourth column shows the corresponding FETI-DPH CPU time in seconds.
    The signs ``-" and ``+" indicate FETI-DPH solves on $\Kbold-\frac{1}{i\sqrt{\phi}}\Vbold$ and $\Kbold+\frac{1}{i\sqrt{\phi}}\Vbold$, respectively.
    The fifth and sixth columns show the number of FETI-DPH iterations with edge-based rigid body modes augmentation,
    and the seventh column shows the corresponding FETI-DPH CPU time in seconds.
    Table~\ref{ta:FETIDPphidependence} shows that, in the absence of augmentation, 
    the number of iterations tends to increase, but not strictly, as $\phi$ decreases.
		This dependence on $\phi$ is alleviated by the introduction of rigid body modes augmentation and some CPU time is saved.
    Table~\ref{ta:FETIDPphidependence} demonstrates that rigid body modes augmentation works well for the thermal optimal control problem.
    The effects of rigid body modes augmentation is even more dramatic in the case of solid elements, as described in the following section.
    Further research on the existence of an optimal augmentation for the thermal optimal control problem is an interesting future topic
    (e.g., a set of modes that decreases the number of FETI-DPH iterations to an order of a constant regardless of the values of $\phi$).
		Table~\ref{ta:FETIDPphidependence} also shows that almost the same number of iterations of FETI-DPH
		is required on $\Kbold-\frac{1}{i\sqrt{\phi}}\Vbold$ and $\Kbold+\frac{1}{i\sqrt{\phi}}\Vbold$ for each value of $\phi$. 
    This result is not surprising, given that the two operators are complex conjugate to each other, and consequently have complex conjugate eigenvalues, 
    the same pattern of clusterings, and an identical condition number.

    There are two kinds of scalability tests for domain decomposition methods: 
    numerical scalability and parallel scalability \cite{Farhat2005}.
    Numerical scalability is studied by varying the problem size, subdomain size, and the number of subdomains.
		Table~\ref{ta:FETIDPnumericalscalability} shows numerical scalability of FETI-DPH on the thermal problem~\eqref{eq:heatconduction}.
    The number of iterations is shown for various mesh sizes $h$ and for various subdomain sizes $H$.
    The regularization parameter $\phi$ is $2\times10^{-8}$.
    The convergence threshold for FETI-DPH is $10^{-10}$.
    The fixed ratio $H/h = 2^6$ is used. 
    Because the ratio $H/h$ is fixed to be $2^6$, each subdomain contains 4096 elements.
    The problem size is varied from around 15 thousand degrees of freedom to 67 million degrees of freedom.
    Based on the theoretical result of \eqref{eq:conditionnumberbound}, the number of FETI-DPH iterations must be more or less
    constant no matter what problem size is considered, provided that the ratio $H/h$ is fixed.
    Indeed, Table~\ref{ta:FETIDPnumericalscalability} shows that it is the case.
		As $h$ decreases, the total number of elements increases, meaning that the problem size increases.
		Table~\ref{ta:FETIDPnumericalscalability} shows that the number of FETI-DPH iterations actually decreases as the problem size increases,
    which is consistent with \eqref{eq:conditionnumberbound}.
		Again, almost the same number of iterations of FETI-DPH 
    is required on $\Kbold-\frac{1}{i\sqrt{\phi}}\Vbold$ and $\Kbold+\frac{1}{i\sqrt{\phi}}\Vbold$ for each $h$.	

    \begin{table}[th]
    \caption{Numerical scalability of the FETI-DPH solver for a fixed $\phi$ value for the heat control problem.}
    \centering
		\footnotesize
    \begin{tabular}{|r|r|c|c|c|c|c|}
    \hline  num. of dofs      &  num. of elem.   &  h          &  H      & nsub & $\Kbold-\frac{1}{i\sqrt{\phi}}\Vbold$ & $\Kbold+\frac{1}{i\sqrt{\phi}}\Vbold$
 \\ \hline
    \hline       15,876       &      16,384      & $1/128$     &  $1/2$  &  4   &          26                           &         25
 \\ \hline       64,516       &      65,536      & $1/256$     &  $1/4$  &  16  &          25                           &         25
 \\ \hline      260,100       &     262,144      & $1/512$     &  $1/8$  &  64  &          21                           &         20
 \\ \hline    1,044,484       &   1,048,576      & $1/1024$    &  $1/16$ &  256 &          15                           &         15
 \\ \hline    4,186,116       &   4,194,304      & $1/2048$    &  $1/32$ & 1024 &          13                           &         13
 \\ \hline   16,760,836       &  16,777,216      & $1/4096$    &  $1/64$ & 4096 &          12                           &         13
 \\ \hline   67,076,100       &  67,108,864      & $1/8192$    & $1/128$ &16384 &          12                           &         12 
 \\ \hline
    \end{tabular}
    \label{ta:FETIDPnumericalscalability}
    \end{table}

    Parallel scalability measures how fast a domain decomposition method converges for a fixed problem size, 
    a fixed subdomain size, and a fixed number of subdomains with increasing number of processes.
    Table~\ref{ta:FETIDPparallelscalability} presents the parallel scalability of FETI-DPH for 
		the heat conduction control problem~\eqref{eq:heatconduction}.
    A mesh size $h=1/1024$ (around 1 million degrees of freedom) and subdomain size $H=1/16$ 
    (4096 elements in each subdomain) are used. 
    The regularization parameter $\phi$ is $2\times10^{-8}$.
    The convergence threshold is $10^{-10}$.
    The number of processes ($N_p$) varies from 1 to 144 
    and parallel speed-up is measured relative to the FETI-DPH CPU time of $N_P=1$.
		As the number of processes increases up to $144$, 
    Table~\ref{ta:FETIDPparallelscalability} shows that a parallel speed-up of $26.8$ is gained. 

    \begin{table}[th]
    \caption{Parallel scalability of the FETI-DPH solver for the heat control problem.}
    \centering
    \begin{tabular}{|r|c|c|}
    \hline         $N_p$      &  CPU time (sec)  &  Parallel speed-up
 \\ \hline
    \hline           1        &    101.99        &       1
 \\ \hline           2        &     52.60        &      1.9
 \\ \hline           4        &     28.13        &      3.6
 \\ \hline           8        &     15.51        &      6.6
 \\ \hline           16       &      9.02        &     11.3 
 \\ \hline           32       &      6.02        &     16.9
 \\ \hline           64       &      5.52        &     18.5
 \\ \hline           128      &      4.30        &     23.7
 \\ \hline           144      &      3.81        &     26.8
 \\ \hline
    \end{tabular}
    \label{ta:FETIDPparallelscalability}
    \end{table}

    The thermal element (i.e., finite element for Laplacian operator) is widely used 
    when one wants to verify a new numerical method or to analyze it.
    In the next section, a more difficult problem is considered (in the sense that properties of a real
    material are used and 3D solid elements are used).

  \subsection{Solid Cantilever Control}
  \label{sec:cantileverproblem}
\begin{figure}[t]
    \centering
    $\begin{array}{c@{\hspace{0.1in}}c}
      \includegraphics[scale=0.17]{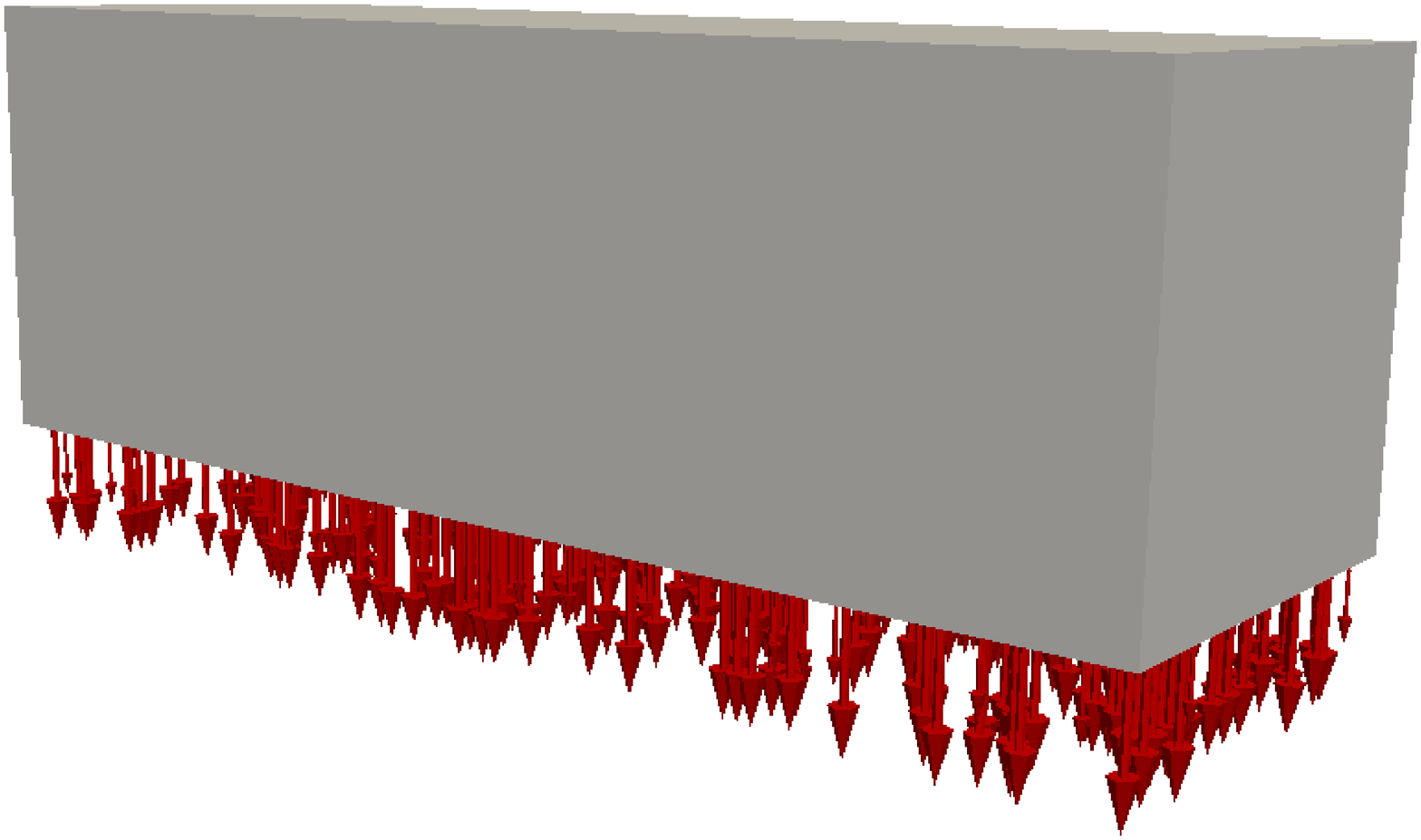} & 
      \includegraphics[scale=0.17]{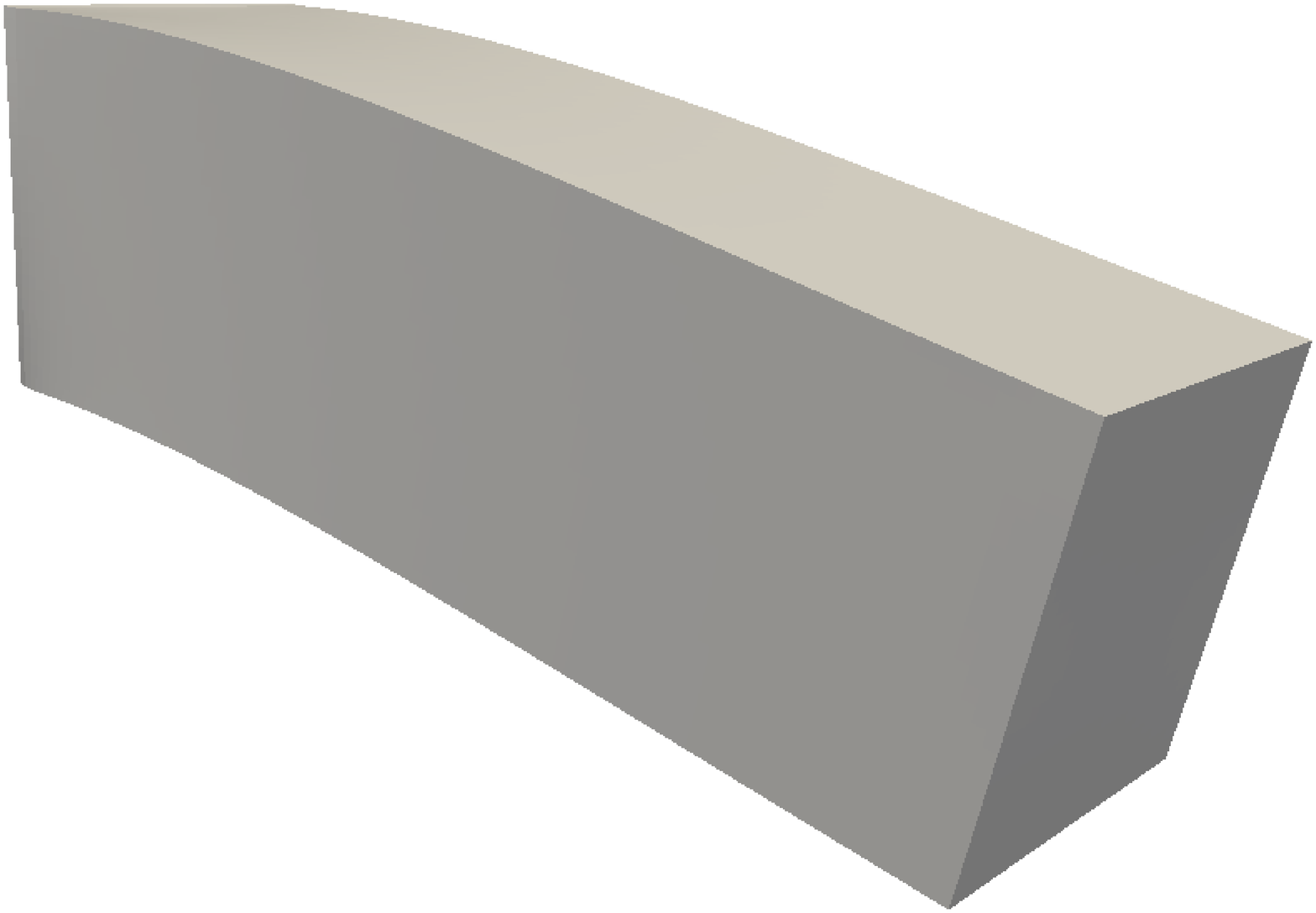} \\
    \mbox{\bf (a) } & \mbox{\bf (b)}
    \end{array}$
    \caption{\textbf{(a)} a uniform downward pressure of 100 kPa is applied to the undeformed configuration of cantilever.
             \textbf{(b)} the deformed configuration that is used as a target state in the optimal control problem.}
    \label{fig:cantilever}
    \end{figure}

    In this section, the proposed optimal control techniques are applied to a linear static structural problem, a cantilever with solid elements.
    In this case, the constraints in \eqref{eq:heatconduction} are replaced by an elastostatic PDE.
    The cantilever has a square cross-section of 1 $\times$ 1 $\text{m}^2$ and 
    a length of 3 m. The target displacement $\bar{\ybold}$ is generated by running 
    a forward PDE simulation with a uniform pressure load of 100 kPa applied to the bottom 
    surface of the cantilever. The following material properties are used:
    Young's modulus of $20.7$ MPa, Poisson's ratio of $0.45$, and density of $1.1 \text{ kg/m}^3$.
    The applied force and the initial configuration of the cantilever are shown on 
    Figure \ref{fig:cantilever}(a) and the corresponding target deformed configuration 
    on Figure~\ref{fig:cantilever}(b). Note that the left end is completely fixed. 
    The body force is used as a control variable. Although this control variable is not practical
    because the body force control is not physically attainable,
    this problem is useful to investigate how the FETI-DPH solver performs on each factor of the Schur complement factorization
    for problems involving solid finite elements and linear elastostatic PDEs.

    \begin{figure}[h]
    \centering
    $\begin{array}{c@{\hspace{0.1in}}c}
      \includegraphics[width=6.4cm, height=4cm]{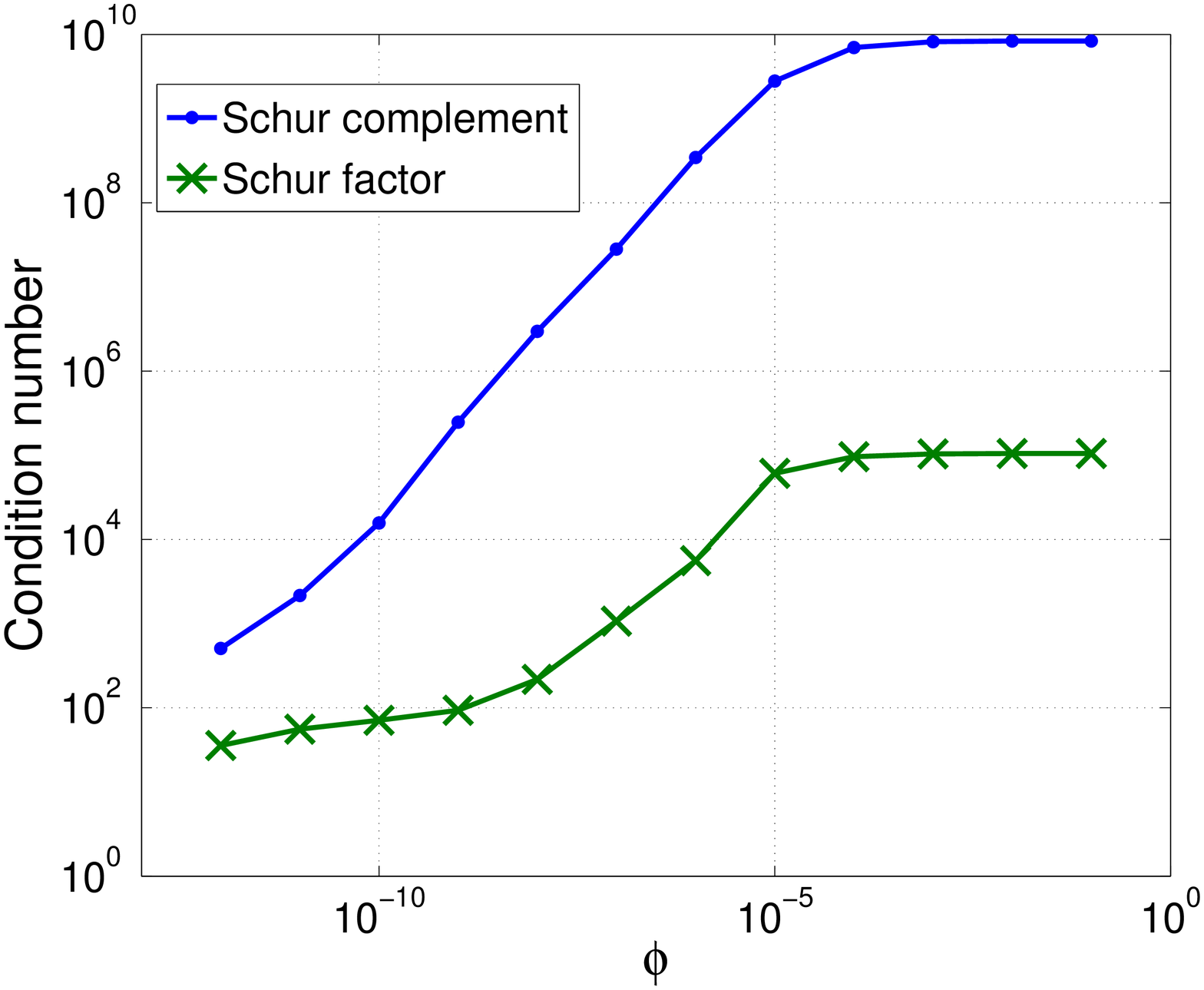} &
      \includegraphics[width=6.4cm, height=4cm]{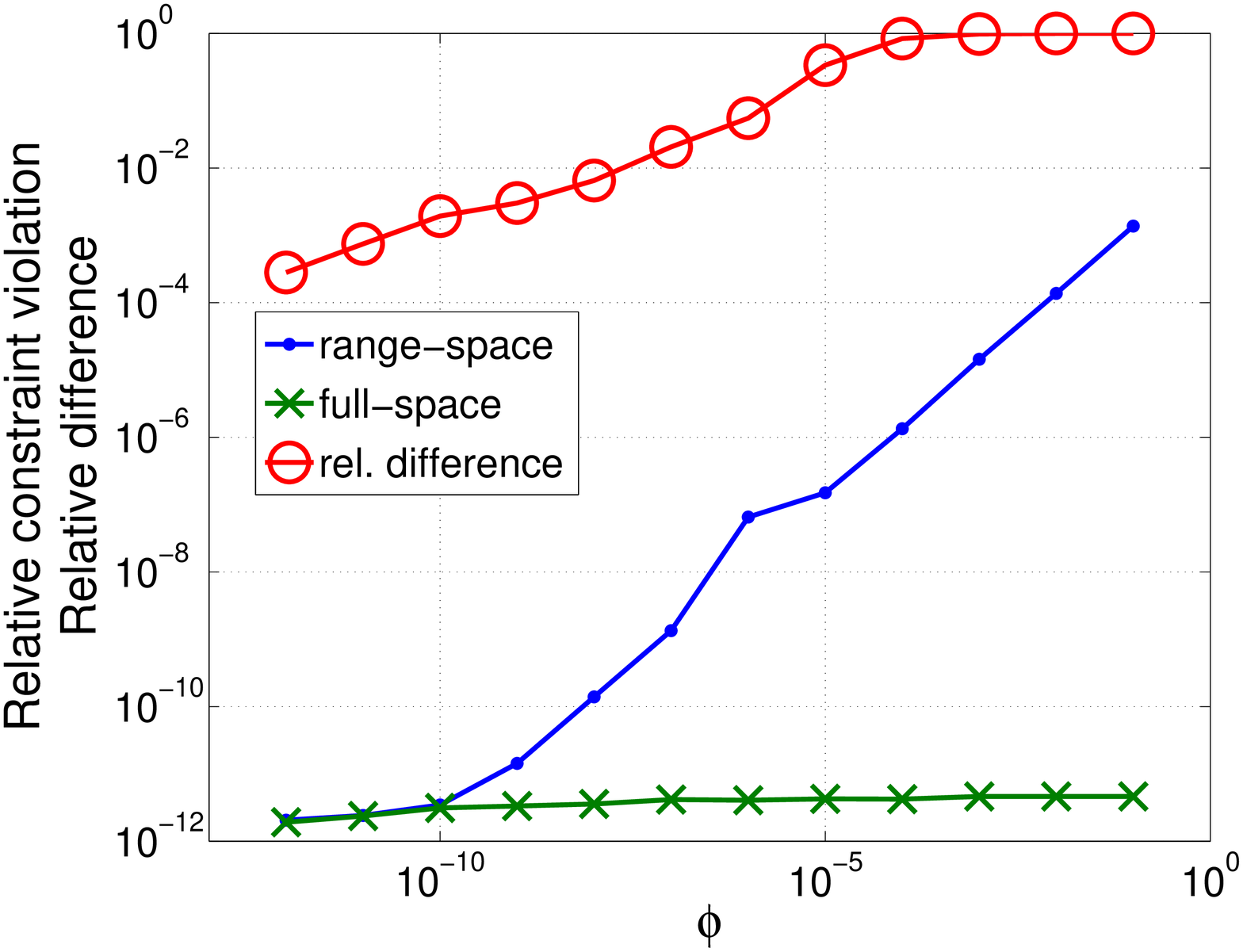} \\
    \mbox{\bf (a) } & \mbox{\bf (b)}
    \end{array}$
    \caption{\textbf{(a)} Condition number vs.\ $\phi$
             \textbf{(b)} Relative constraint violation and relative difference between the target and solution displacements.}
    \label{fig:rangespaceErrorAnalysis}
    \end{figure}

    Figure \ref{fig:rangespaceErrorAnalysis} shows the accuracy issue of the range-space method that was explained in 
    Section \ref{subsec:rangespacemethod}. The results were generated for a relatively small problem ($h = 1/8$ and 5184 degrees of freedom). 
    Figure~\ref{fig:rangespaceErrorAnalysis}(a) shows how the condition numbers of both the negative 
    Schur complement $\Sbold=\Kbold\Vbold^{-1}\Kbold+\frac{1}{\phi}\Vbold$ and a complex factor $\Kbold+\frac{1}{i\sqrt{\phi}}\Vbold$ vary as $\phi$ varies.
    Both condition numbers increase as $\phi$ increases but the order of magnitude of the condition number of the factor
    is one half that of $\Sbold$. Figure~\ref{fig:rangespaceErrorAnalysis}(b) demonstrates the accuracy dependency
    of both the range-space and the full-space methods on the value of $\phi$ by showing the relative constraint violation
    marked by a blue line with dots and a green line with x symbol 
    (see the constraint in \eqref{eq:discretizedoptimalheatcontrol}). The relative constraint violation is defined as 
    \begin{equation}
    \label{eq:relativeconstraintviolation}
       \frac{\|\Kbold\ybold + \Kbold_{c}\ybold_c - \Vbold\ubold\|_2}{\|\Kbold\ybold\|_2}.
    \end{equation}
    As $\phi$ increases, the accuracy of the range-space method degenerates. 
    For example, for values of $\phi$ larger than $10^{-2}$, the constraint violation becomes larger than $10^{-4}$,
    which is of questionable acceptability. On the other hand, the full-space method shows
    high accuracy consistently for the entire range of $\phi$ values considered.
    It is possible to achieve small improvements in the accuracy of the range-space method 
    by using iterative refinement and higher floating-point precision. 
    However, this approach incurs an additional cost and preliminary experiments suggest that it is not a competitive solution. 
    One could argue that the values of interest of the parameter $\phi$ 
    for this problem are the ones smaller than $10^{-5}$ because the relative difference between
    target and solution states is greater than 83$\%$ for values of $\phi$ outside this range. For this range of $\phi$ values,
    the accuracy of the range-space method is acceptable for this particular mesh size ($h=1/8$).

    \begin{table}[h]
    \caption{FETI-DPH's $\phi$ dependency for the structure control problem.}
    \centering
    \footnotesize
    \begin{tabular}{|c||c|c|c||c|c|c||c|}
         \hline   $\phi$ &  -  &  +  &  CPU   & RAUG-& RAUG+  & CPU    & A. speed-up
      \\ \hline
         \hline   1E-7   & 327 & 344 & 909.9  & 30   & 33     & 187.4  &  4.9
      \\ \hline   1E-8   & 341 & 357 & 880.1  & 29   & 32     & 189.4  &  4.6 
      \\ \hline   1E-9   & 286 & 300 & 777.4  & 29   & 31     & 192.4  &  4.0
      \\ \hline   1E-10  & 244 & 258 & 690.2  & 29   & 31     & 176.1  &  3.9
      \\ \hline   1E-11  & 228 & 239 & 673.9  & 29   & 30     & 185.3  &  3.6
      \\ \hline   1E-12  & 225 & 235 & 657.2  & 28   & 29     & 180.1  &  3.6
      \\ \hline   1E-13  & 220 & 229 & 652.3  & 27   & 28     & 193.6  &  3.4
      \\ \hline   1E-14  & 236 & 214 & 682.9  & 29   & 30     & 186.6  &  3.7
      \\ \hline   1E-15  & 262 & 260 & 773.2  & 34   & 34     & 173.5  &  4.5
      \\ \hline   1E-16  & 280 & 284 & 817.6  & 40   & 39     & 222.3  &  3.7
      \\ \hline   1E-17  & 315 & 316 & 976.3  & 46   & 46     & 214.1  &  4.6
      \\ \hline   1E-18  & 321 & 321 & 1009.6 & 52   & 52     & 236.5  &  4.3
      \\ \hline
    \end{tabular}
    \label{ta:FETIDPphidependence_structure}
    \end{table}

    Table~\ref{ta:FETIDPphidependence_structure} shows FETI-DPH's dependence on the value of $\phi$ for the cantilever control problem.
    The mesh size is $h=1/90$ and the convergence threshold is $10^{-9}$.
    The size of the subdomain is $H=1/6$ and the number of processes is 64.
    The second and third columns show the number of FETI-DPH iterations without any augmentation, and
    the fourth column shows the corresponding FETI-DPH CPU time in seconds.
    The signs ``-" and ``+" indicate the FETI-DPH solves on $\Kbold-\frac{1}{i\sqrt{\phi}}\Vbold$ and $\Kbold+\frac{1}{i\sqrt{\phi}}\Vbold$, respectively.
    The fifth and sixth columns show the number of FETI-DPH iterations with edge-based rigid body modes augmentation, 
    and the seventh column shows the corresponding FETI-DPH CPU time in seconds.
    The last column shows the speedup due to augmentation.
    Without augmentation, the performance of the FETI-DPH solver degrades as the value of $\phi$ increases from 
    $10^{-12}$ to $10^{-7}$. This makes sense because the condition number  
    of the complex factors in \eqref{eq:factorizationofschur} and the Schur complement itself increase as $\phi$ increases
    (see Figure \ref{fig:rangespaceErrorAnalysis}). However, for the range of smaller $\phi$ values (i.e., less than $10^{-12}$),
    the number of FETI-DPH iterations without augmentation increases as $\phi$ decreases. 
    This is analogous to FETI-DPH being dependent on wave numbers in acoustic problems if no augmentation is applied
    (i.e., more iterations are required for a larger wave number without augmentation) although the complex factors 
    are not the same as in the system that normally arises from a Helmholtz problem.
    In order to alleviate this deterioration, edge-based rigid body mode augmentation is used. 
    Table~\ref{ta:FETIDPphidependence_structure} shows that the augmentation decreases the number of FETI-DPH iterations
    substantially and reduces the CPU time by a factor of four on average.
    However, the rigid body modes augmentation is not optimal in the sense that
    as $\phi$ decreases the number of FETIP-DPH iterations still increases.

    Table \ref{ta:FETIscalability_structure} shows the numerical scalability of the FETI-DPH solver for the cantilever control problem.
		The number of iterations is shown for various mesh sizes $h$ and for various subdomain sizes $H$.
    The regularization parameter $\phi$ is $1\times10^{-16}$.
    The convergence threshold for FETI-DPH is $10^{-9}$.
    The problem size is varied from around 2 million degrees of freedom to 40 million degrees of freedom.
    As in the numerical scalability test for the thermal problem, the $H/h$ ratio is fixed to 15 (exactly 3375 elements in each subdomain) 
    and the number of FETI-DPH iterations is counted for both factors $\Kbold-\frac{1}{i\sqrt{\phi}}\Vbold$ and $\Kbold+\frac{1}{i\sqrt{\phi}}\Vbold$. 
    As the problem size increases, the number of FETI iterations required for convergence slightly decreases, which is
    again consistent with the theoretical result~\eqref{eq:conditionnumberbound}.

    \begin{table}[h]
    \caption{Numerical scalability of the FETI-DPH solver for a fixed $\phi$ value for the structure control problem.}
    \centering
    \footnotesize
    \begin{tabular}{|r|r|c|c|c|c|c|}
    \hline  num.\ of dofs  &  num.\ of elem.   &  h          &  H      & nsub & $\Kbold-\frac{1}{i\sqrt{\phi}}\Vbold$ & $\Kbold+\frac{1}{i\sqrt{\phi}}\Vbold$
 \\ \hline
    \hline    1,976,400   &       648,000    & $1/60$      &  $1/4$  &  192 &          43                 &         43
 \\ \hline    3,847,500   &     1,265,625    & $1/75$      &  $1/5$  &  375 &          41                 &         40    
 \\ \hline    6,633,900   &     2,187,000    & $1/90$      &  $1/6$  &  648 &          39                 &         38
 \\ \hline   10,517,850   &     3,472,875    & $1/105$     &  $1/7$  & 1029 &          37                 &         37
 \\ \hline   15,681,600   &     5,184,000    & $1/120$     &  $1/8$  & 1536 &          36                 &         36    
 \\ \hline   22,307,400   &     7,381,125    & $1/135$     &  $1/9$  & 2187 &          35                 &         34 
 \\ \hline   30,577,500   &    10,125,000    & $1/150$     &  $1/10$ & 3000 &          33                 &         33      
 \\ \hline   40,674,150   &    13,476,375    & $1/165$     &  $1/11$ & 3993 &          33                 &         33 
 \\ \hline
    \end{tabular}
    \label{ta:FETIscalability_structure}
    \end{table}

    Table~\ref{ta:FETIDPparallelscalability_structure} shows the parallel scalability of FETI-DPH for the cantilever control problem. 
    The FETI-DPH CPU time and the corresponding parallel speed-up are shown for the various numbers of processes.
    The regularization parameter $\phi$ is $10^{-12}$.
    The convergence threshold is $10^{-9}$.
    A mesh size $h=1/90$ (i.e., around 6.6 million degrees of freedom) and a subdomain size $H=1/6$ are used.
    The number of processes ($N_p$) increases from 1 to 64
    and parallel speed-up is measured as the ratio with respect to the CPU time of $N_P=1$.
    As the number of processes increases up to $64$,
    Table~\ref{ta:FETIDPparallelscalability_structure} shows the speed up is greater than half the maximum possible.

    \begin{table}[h]
    \caption{Parallel scalability of FETI-DPH solver for the structure control problem.}
    \centering
    \begin{tabular}{|r|c|c|}
    \hline         $N_p$      &  CPU time (sec)  &  Parallel speed-up
 \\ \hline
    \hline           1        &    5854.9        &     1.0  
 \\ \hline           2        &    2834.6        &     2.1 
 \\ \hline           4        &    1547.0        &     3.8 
 \\ \hline           8        &     886.6        &     6.6 
 \\ \hline           16       &     586.4        &    10.0 
 \\ \hline           32       &     317.4        &    18.4 
 \\ \hline           64       &     163.5        &    35.8
 \\ \hline
    \end{tabular}
    \label{ta:FETIDPparallelscalability_structure}
    \end{table}

    Due to the introduction of complex numbers in the factorization, more storage and more computational cost per solve 
    are required compared to Pearson and Wathen's approximation to the Schur complement ($\Sbold_p$ in \eqref{eq:pearsonapproximation}).
    Thus, one must consider the effectiveness of the complex factorization carefully although only one
    solve with complex factorization is required in the range-space method and fewer iterations are
    required in the full-space method. Table~\ref{ta:Sc_Sp_comparisonFETI} reports 
    comparison of CPU time and number of iterations between $\Sbold_p$ and $\Sbold_c$ in the full-space method.
    The block diagonal preconditioner of Murphy, et al.\ \eqref{eq:murphypreconditioner} is used.
    Comparison is made for various values of $\phi$ and each solve is done by FETI-DPH with edge augmentation.
    The mesh size is $h=1/30$ and the convergence threshold for FETI-DPH is $10^{-9}$.
    GMRES is used for the main solver and the convergence threshold for GMRES is $10^{-10}$.
    The size of each subdomain is $H=1/2$ and 64 processes are used.
    The second column ($B_{\Sbold_p}$) shows the CPU time for building an operator and
    the third column ($S_{\Sbold_p}$) shows the CPU time for one solve with $\Sbold_p$.
    The fourth column ($NI_{\Sbold_p}$) shows the number of solves with $\Sbold_p$ that are required for convergence.
    The fifth column (TCPU) shows the total CPU time in seconds.
    The sixth to ninth columns show the corresponding results for $\Sbold_c$.
    Note that the number of GMRES iterations required for convergence with $\Sbold_c$ is more than 3, which is not consistent with
    the spectral analysis done in \cite{Murphy2000}. 
    Each solve time and building time is considerably higher for $\Sbold_c$ than $\Sbold_p$, as expected. 
    However, the number of GMRES iterations
    required for convergence is much higher for $\Sbold_p$ than $\Sbold_c$. Additionally, the number of iterations for $\Sbold_p$ increases
    as $\phi$ decreases, while for $\Sbold_c$ it is bounded above. 
    For this particular problem $\Sbold_p$ is a better choice for relatively large values of $\phi$ 
    (i.e., $\phi=10^{-7}$), while $\Sbold_c$ is a better choice for any value of $\phi$ smaller than $10^{-9}$.

    \begin{table}[h]
    \caption{Comparison between $\Sbold_p$ and $\Sbold_c$ in the full-space method.}
    \centering
    \footnotesize
    \begin{tabular}{|c||c|c|c|c||c|c|c|c|}
         \hline   $\phi$ & $B_{\Sbold_p}$ & $S_{\Sbold_p}$ & $NI_{\Sbold_p}$ & TCPU   & $B_{\Sbold_c}$ & $S_{\Sbold_c}$ & $NI_{\Sbold_c}$ & TCPU
      \\ \hline
         \hline   1E-7   &   2.1     &  2.0      &     43     & 88.1   &   9.2     &   3.9     &    25      & 106.5
      \\ \hline   1E-8   &   2.1     &  1.9      &     49     & 95.2   &   8.9     &   3.7     &    25      & 101.4
      \\ \hline   1E-9   &   2.1     &  2.0      &     55     & 112.1  &   9.6     &   4.1     &    25      & 112.1
      \\ \hline   1E-10  &   2.3     &  2.3      &     75     & 174.8  &   9.6     &   3.9     &    25      & 107.1
      \\ \hline   1E-11  &   2.4     &  2.4      &     77     & 187.2  &   9.3     &   4.1     &    25      & 111.8
      \\ \hline   1E-12  &   2.1     &  2.0      &     97     & 196.1  &   8.7     &   3.9     &    25      & 106.2
      \\ \hline   1E-13  &   2.2     &  2.0      &     97     & 196.2  &   9.2     &   4.6     &    25      & 124.2
      \\ \hline   1E-14  &   2.2     &  2.1      &     97     & 205.9  &   9.6     &   6.0     &    23      & 147.6
      \\ \hline   1E-15  &   2.2     &  2.1      &    110     & 233.2  &   8.8     &   5.6     &    23      & 137.6
      \\ \hline   1E-16  &   2.3     &  2.3      &    323     & 745.2  &   9.0     &   6.9     &    23      & 167.7
      \\ \hline
    \end{tabular}
    \label{ta:Sc_Sp_comparisonFETI}
    \end{table}

	\section{Conclusion}
	\label{sec:conclusion}
    We have introduced a practical factorization of the Schur complement that arises from distributed optimal control of linear static systems.
    Due to the exact representation, if the range-space method is applicable, then one solve with the Schur complement is
    sufficient to obtain an optimal control solution. However, the Schur complement becomes ill-conditioned 
    for large values of $\phi$ and as the mesh is refined. For example, the Schur complement that arises
    from a 3D cantilever problem with real material properties of rubber is prone to ill-conditioning when a relatively large regularization value is used. 
    In such a case, an inaccurate solution is likely to be obtained if the range-space method is used, 
    so solving a full KKT system simultaneously using a Krylov iterative method with a good preconditioner is recommended.
    The complex factorization of the Schur complement introduced in this paper can be used with a Schur-complement based preconditioner.
    The comparison between the approximate Schur complement of Pearson and Wathen and the complex factorization as a preconditioner
    in the full-space method shows promising results for the complex factorization even in the context of the full-space method.
    The Schur complement solve is done with the parallel domain decomposition linear solver FETI-DPH.   
    The scalability of FETI-DPH (both numerical and parallel) as well as its dependence on $\phi$
    are studied in two academic problems: a thermal 2D problem and a structural 3D problem.  
    Edge-based rigid body modes augmentation is able to bring the number of iterations down, but further research is necessary  
    to find an optimal augmentation in the context of PDE-constrained distributed optimal control.
    The combination of exact representation of the Schur complement and good scalability of the FETI-DPH solver
    in addition to extensibility of the representation indicates promise for use of the complex factorization
    in more complicated and practical problems.


	\section*{Acknowledgments}
		The authors thank Philip Avery in the Farhat Research Group for his valuable comments and 
    essential help with coding the physics-based C++ PDE solver Aero-S.

	\bibliographystyle{plain}
	\bibliography{citations}

\end{document}